\documentclass[11pt]{amsart}
\usepackage{geometry}                % See geometry.pdf to learn the layout options. There are lots.
\usepackage{graphicx}
\usepackage{amssymb}
\usepackage{epstopdf}
\usepackage{hyperref}
\usepackage{cleveref}
\usepackage{todonotes}
\usepackage{caption,subcaption}
\usepackage{enumerate}
\usepackage{amscd}
\usepackage{accents}

\makeatletter
\@namedef{subjclassname@2020}{%
  \textup{2020} Mathematics Subject Classification}
\makeatother

\def\endpf{\relax\ifmmode\expandafter\endproofmath\else
  \unskip\nobreak\hfil\penalty50\hskip.75em\hbox{}\nobreak\hfil\bull
  {\parfillskip=0pt \finalhyphendemerits=0 \bigbreak}\fi}
\def\bull{\vbox{\hrule\hbox{\vrule\kern3pt\vbox{\kern6pt}\kern3pt\vrule}\hrule}}

\DeclareMathOperator{\Int}{int}

\newtheorem{defn}{Definition}[section]
\newtheorem{lemma}[defn]{Lemma}

\newtheorem{theorem}[defn]{Theorem}
\newtheorem{definition}[defn]{Definition}

\newtheorem{proposition}[defn]{Proposition}
\newtheorem{corollary}[defn]{Corollary}

\newtheorem{maintheorem}{Theorem}

\newtheorem{example}[defn]{Example}

\newcommand{\cc}{{\mathbb C}}
\newcommand{\zz}{{\mathbb Z}}
\newcommand{\rr}{{\mathbb R}}

\newcommand{\pp}{{\mathbb P}}

\newcommand{\spin}{\ifmmode{\rm Spin}\else{${\rm spin}$\ }\fi}
\newcommand{\spinc}{\ifmmode{{\rm Spin}^c}\else{${\rm spin}^c$\ }\fi}

\newcommand{\DC}{\mathcal{DC}}
\renewcommand{\L}{\mathcal{L}}

\newcommand{\calh}{\mathcal{H}}

\newcommand{\calf}{\mathcal{F}}
\newcommand{\calV}{\mathcal{V}}

\def\veps{\varepsilon}
\def\Rmv{{R - \varepsilon}}
\def\Rpv{{R + \varepsilon}}
\def\vfi{\varphi}

\newif\ifpic
\picfalse    % no pictures
\pictrue   % with pictures

\DeclareMathOperator{\lk}{lk}
\DeclareMathOperator{\sign}{sign}

\begin{document}

\title{Disoriented homology and double branched covers}
\author[Brendan Owens]{Brendan Owens} 
\address{School of Mathematics and Statistics, University of Glasgow, University Place, Glasgow G12 8QQ, United Kingdom}
\email{brendan.owens@glasgow.ac.uk}

\author[Sa\v{s}o Strle]{Sa\v{s}o Strle}
\address{Faculty of Mathematics and Physics, University of Ljubljana, Jadranska 19, 1000 Ljubljana, Slovenia}
\email{saso.strle@fmf.uni-lj.si}

\date{\today}                                           % Activate to display a given date or no date

%; Secondary 57R19}

\begin{abstract}
This paper provides a convenient and practical method to compute the homology and intersection pairing of a  branched double cover of the 4-ball.

To projections of links in the 3-ball, and to projections of surfaces in the 4-ball into the boundary sphere, we associate a sequence of homology groups, called the disoriented homology.  We show that the disoriented homology is isomorphic to the homology of the double branched cover of the link or surface.  We  define a pairing on the first disoriented homology group of a 
surface and show that this is equal to the intersection pairing of the branched cover.  These results generalize work of Gordon and Litherland, for embedded surfaces in the 3-sphere, to arbitrary surfaces in the 4-ball.  We also give a generalization of the signature formula of Gordon-Litherland to the general setting.

Our results are underpinned by a theorem describing a handle decomposition of the branched double cover of a codimension-2 submanifold in the $n$-ball, which generalizes previous results of Akbulut-Kirby and others.
\end{abstract}

\maketitle

\section{Introduction}

Branched covering spaces have proved to be an extremely efficient way of encoding embedding information about submanifolds \cite{ak,greene,lisca,OSz}. The basic information about a covering space is its homology; this is often the starting point for extracting other invariants, such as various gauge theoretic invariants. In \cite{gl}, Gordon and Litherland showed that the first homology of an embedded spanning surface $F$ for a link $L$ in $S^3$ is isomorphic to the second homology of the double cover $X$ of the 4-ball branched along the properly-embedded surface obtained by pushing the interior of $F$ into the ball. Moreover, they defined a bilinear form on $H_1(F)$ and showed that it is isomorphic to the intersection form of $X$; they also derived a formula for the signature of $L$ in terms of this form. 

The main goal of this paper is to generalize these results to embedded surfaces in the 4-ball. As a warm-up we consider links and tangles in the 3-ball.  We use the radial distance function to induce a bridge decomposition on the radial projection $P \subset S^2$ of the link or tangle $L \subset B^3$, as in the example of the trefoil shown in 
\begin{figure}[htbp]
\centering
\includegraphics[scale=1]{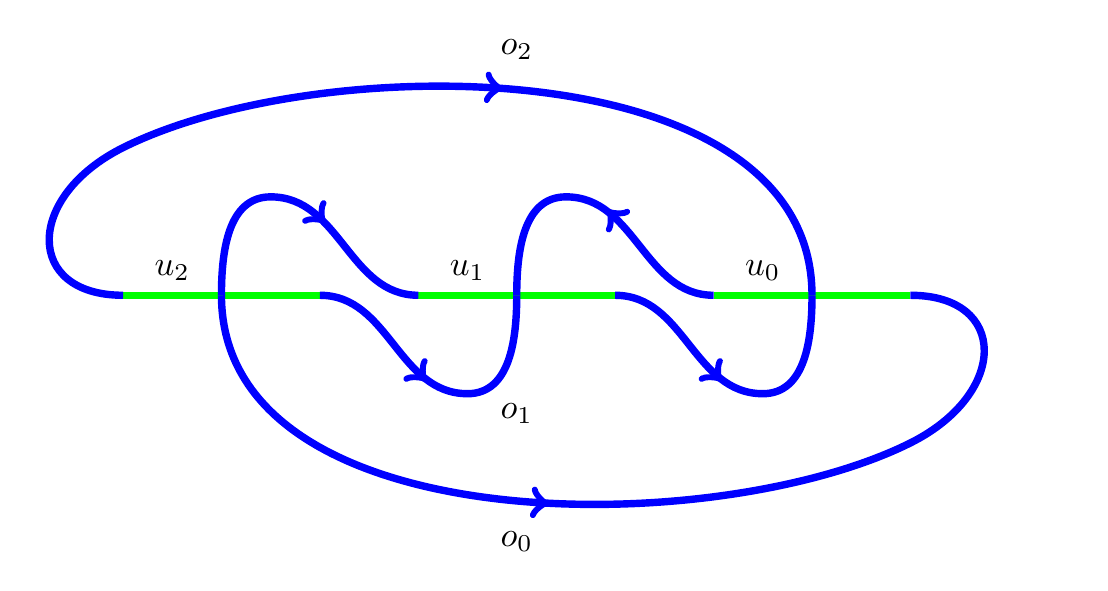}
\caption
{{\bf A bridge decomposition of the left-handed trefoil.} Underbridges are shown in green, with overbridges in blue.  The arrows on the overbridges specify a disorientation.}
\label{fig:trefoilDH}
\end{figure}
Figure \ref{fig:trefoilDH}.  We choose a disorientation of each overbridge, again as shown in the figure: each segment of the complement of the underbridges in the overbridge is given an orientation.  These are chosen so that the orientation switches at each crossing.  

We use this data to define the disoriented chain complex $\DC_*(P)$ (see \Cref{sec:tangle}): $\DC_1(P)$ is the free abelian group generated by the overbridges,  $\DC_0(P)$ is generated by the underbridges, and the boundary operator between them is given by counting with sign how many times each overbridge points into or out of each underbridge.  The boundary operator from $\DC_0(P)$ to $\DC_{-1}(P)=\zz$ is the augmentation homomorphism. We show that the homology of this complex computes that of the double branched cover of $L$ (for a precise statement see \Cref{prop:3d-homology}; this is closely related to the fact that the \emph{coloring matrix} of a link diagram presents the first homology of the double branched cover, cf. \cite{josh,liv}):

\begin{maintheorem}
\label{thm:tangleintro}
The disoriented homology of a link or tangle $L$ in $B^3$ is isomorphic to the shifted reduced homology of the double cover of $B^3$ branched along $L$, i.e.,
$$H_*(\DC_*(P)) \cong \widetilde H_{*+1}(\Sigma_2(B^3,L)).$$
\end{maintheorem}

\begin{figure}[htbp]
\centering
\includegraphics[scale=1.5]{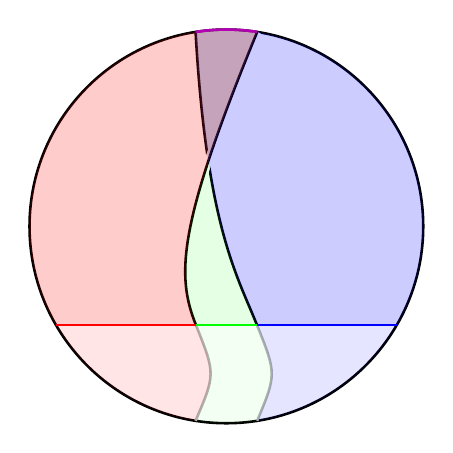}
\caption
{{\bf The radial projection of a projective plane with a compatible handle decomposition.} The round disk is the 0-handle, the green band the 1-handle, and the red and blue disks combine to give the 2-handle. The green arc signifies the ribbon singularity. The 2-handle is split into four subdisks by its intersection with the ribbon surface.}
\label{fig:proj-plane-intro}
\end{figure}

\begin{figure}[htbp]
\centering
        \includegraphics[scale=1.2]{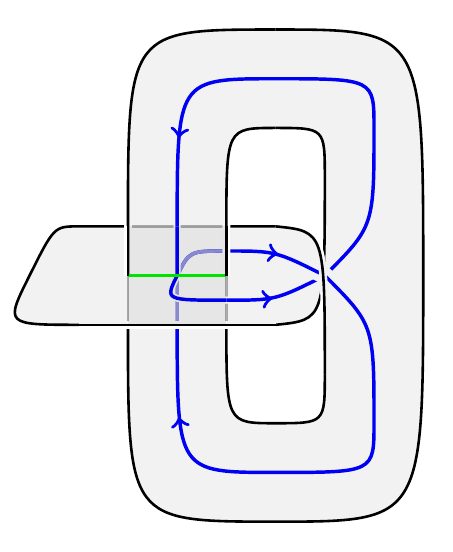}  % [scale=1.2]
\caption
{{\bf A generating set for disoriented homology.} A ribbon-immersed annulus with a generator of its first disoriented homology group. The ribbon singularity is shown in green.}
\label{fig:DHgen}
\end{figure}

The disoriented homology of a compact surface $F$ properly embedded in the 4-ball, with or without boundary, may be defined in a similar manner (see \Cref{sec:slice-disoriented}).  Starting with a handle decomposition of $F \subset B^4$ induced by the radial distance function, we consider the images of these handles in the radial projection $F_s \subset S^3$ of $F$ as handles of $F_s$. Assuming the projection to be regular, $F_s$ may be decomposed as a ribbon-immersed surface and a union of disjoint disks that are 2-handles  of $F_s$, as in \Cref{fig:proj-plane-intro}.
The group $\DC_k(F_s)$ for $k \ge 0$ is freely generated by the $k$-handles of $F_s$.  The boundary operator from $\DC_2$ to $\DC_1$ for each 2-handle essentially counts with sign how many times the intersection of the 2-handle with a 1-handle goes over that 1-handle; this is computed using \emph{disorientations} of handles.   Disorientations of 1-handles are orientations of their cores, switching each time they pass through a ribbon singularity; disorientations of 2-handles are determined by chessboard coloring of the regions into which 2-handles are split by their intersections with the ribbon surface.  The remaining boundary homomorphisms are defined similarly to the 3-dimensional case. We also show that taking linking numbers with double normal pushoffs gives rise to a pairing $\lambda$ on the first disoriented homology group of $F_s$, which we call the GL-pairing of $F_s$.  To define the pairing we use a more geometric description of the first disoriented homology group for a ribbon-immersed surface (see \Cref{sec:ribbon-disoriented}); \Cref{fig:DHgen} shows an example.  We prove the following  (see \Cref{thm:4d-homology} for a more precise formulation):

\begin{maintheorem}
\label{thm:surfaceintro}
The disoriented homology of a properly-embedded compact surface $F$ in the 4-ball is isomorphic to the shifted reduced homology of the  double cover $\Sigma_2(B^4,F)$ of the 4-ball branched along $F$, i.e.,
$$H_*(\DC_*(F_s)) \cong \widetilde H_{*+1}(\Sigma_2(B^4,F)).$$
Moreover, the intersection pairing of $\Sigma_2(B^4,F)$ under this identification agrees with the GL-pairing $\lambda$.
\end{maintheorem}

The proof of this theorem relies on a Kirby diagram for the branched double cover $\Sigma_2(B^4,F)$; in particular, we give a recipe for drawing the attaching spheres of 3-handles. This is illustrated in \Cref{ex:proj-Kirby} for the surface in \Cref{fig:proj-plane-intro}.
%Figure \ref{fig:DHgen} shows a ribbon-immersed surface together with a generator of its first disoriented homology group.

A widely-used application of the celebrated paper of Gordon and Litherland is a convenient formula to compute the signature of a link using the signature of the bilinear form associated to a spanning surface for the link in $S^3$.  We generalise this to give a signature formula based on the GL-pairing of a slice surface as follows:

\begin{maintheorem}
\label{thm:sigintro}
Let a link $\L$ be the boundary of a slice surface $F \subset B^4$, and let $\lambda$ be the GL-pairing on the first disoriented homology group of $F$. Then for any choice $\vec\L$ of orientation for $\L$, its signature is given by
$$\sigma(\vec\L)=\sigma(\lambda) - \frac 12 \lk(\vec\L,\vec\L^F),$$
where $\vec\L^F$ is a parallel copy of $\vec\L$ on the radial projection $F_s$ of $F$, oriented consistently with $\vec\L$.
\end{maintheorem}

The organization of the paper is as follows: we define the disoriented homology of a properly embedded tangle $L \subset B^3$ in \Cref{sec:tangle}.  For a properly embedded surface $F \subset B^4$ we define in \Cref{sec:surface} its description $F_s \subset S^3$ which in the case of a ribbon surface is its ribbon immersion; for a general surface it is its regular projection with a compatible handle decomposition. Based on this description we define the disoriented homology $DH_*(F_s)$, for ribbon surfaces in \Cref{sec:ribbon-disoriented} and for general slice surfaces in \Cref{sec:slice-disoriented}. In case of a ribbon surface the first disoriented homology $DH_1(F_s)$ is a subgroup of the singular homology group $H_1(F_s;\zz)$ of the ribbon-immersed surface $F_s \subset S^3$ generated by those cycles that in a neighborhood of every ribbon sigularity are multiples of the chain pictured in \Cref{fig:localdis}; the structure of this chain also gives the homology its name. We then extend the definition of disoriented homology to general slice surface descriptions $F_s$ and give several alternative descriptions of the groups; see \Cref{def:DHv2} for the description used to relate the disoriented homology of the surface to the homology of the double branched cover as in \Cref{thm:surfaceintro}.

We define the pairing $\lambda$ on $DH_1(F_s)$ in \Cref{sec:pairing}, generalizing the Gordon-Litherland pairing.

\Cref{sec:dbc} is the technical core of the paper, in which we relate a handle decomposition of a codimension-2 submanifold $F$ on the $n$-ball to a handle decomposition of its double branched cover $X$.
In \Cref{sec:dbc3} we show how a bridge decomposition of $L$ gives rise to a handle decomposition of the double cover $Y$ of $B^3$ branched along $L$ and give a recipe for drawing a Heegaard diagram of the double branched cover of a link in $S^3$. We show in \Cref{prop:3d-homology} that the disoriented homology of a tangle is isomorphic to the shifted homology of $Y$, proving \Cref{thm:tangleintro}.
In \Cref{sec:dbc4} we consider the case of a surface $F$ in the 4-ball, and show how to construct a Kirby diagram for $X$ based on a handle decomposition of $F$.  We use this to prove \Cref{thm:surfaceintro}.

In \Cref{sec:signature} we prove \Cref{thm:sigintro}.

\noindent {\bf Acknowledgements:} The authors thank Josh Greene for many helpful conversations.  We thank Frank Swenton for help with TikZ.  We thank the anonymous referee for helpful suggestions to improve the exposition.  The second author was partially supported by Slovenian Research Agency (ARRS) Research program P1-0288.

%%%%%%%%%%%%%%%
\section{Disoriented homology of tangles}\label{sec:tangle}
Let $L$ be a properly embedded compact 1-manifold in the 3-ball, i.e., a tangle or a link, to which the radial distance function $\rho$ restricts to be Morse, giving a handle decomposition of $L$.  This is known as a bridge decomposition of $L$. We assume that the radial projection $P \subset S^2$ of $L$ has only ordinary double points. The bridge decomposition of $L$ induces a bridge decomposition of $P$ which then carries the same information as a diagram of $L$; we refer to double points of $P$ as crossings. In this context $0$-handles and $1$-handles are called underbridges and overbridges respectively.   We further assume that 
\begin{itemize}
\item all endpoints of $P$ are contained in underbridges, and
\item at each crossing, an overbridge crosses over an underbridge.
\end{itemize}

For each overbridge of $P$ choose a \emph{disorientation} as follows: split the overbridge into subarcs separated by crossings and give consecutive subarcs opposite  orientations. Denote the projection with this extra information (bridge decomposition and disorientations of overbridges) by $P^\flat$.  Define the \emph{disoriented chain complex} $\DC_*(P^\flat)$ of $L$ as follows.  Let $\DC_0(P^\flat)$ be the free abelian group generated by the underbridges and
$\DC_1(P^\flat)$ be the free abelian group generated by the disoriented overbridges. The boundary homomorphism $\partial_P^\flat : \DC_1(P^\flat) \to \DC_0(P^\flat)$ associates to each overbridge a linear combination of underbridges, where if an oriented arc of the overbridge points to/from an underbridge, it contributes plus/minus that underbridge.  
Note that the contribution at each crossing is $\pm 2$ times the underbridge at the crossing.
Let $\DC_{-1}(P^\flat) =\zz$ and $\veps : \DC_0(P^\flat) \to \zz$ be the augmentation homomorphism mapping every underbridge to $1$. 
Then 
$$ 0 \to \DC_1(P^\flat) \stackrel{\partial_P^\flat}{\longrightarrow} \DC_0(P^\flat) \stackrel{\epsilon}{\to} \DC_{-1}(P^\flat) \to 0$$
is a chain complex that we refer to as the \emph{disoriented chain complex} of $L$. The homology of this complex is the \emph{disoriented homology} of $L$. We will show in \Cref{sec:dbc3} that this is isomorphic to the shifted homology of the double branched cover of $B^3$ with branch set $L$. In particular this implies that the disoriented homology of $L$ is independent of the choices involved in its definition.

\begin{example}
\label{eg:trefoilDH}
We illustrate the above with the example of the trefoil $L$ as in Figure \ref{fig:trefoilDH}, using the bridge decomposition and the disorientations of the overbridges as in that figure.
Relative to the labellings of the underbridges and overbridges, the boundary homomorphism is given by the matrix
$$\partial_P^\flat=\left[  
\begin{matrix} 1 & -1 & 2 \\ 1 & 2 & -1 \\ -2 & -1 & -1 \end{matrix}
\right].$$
Since the rank of this matrix is $2$, it follows that $H_1(\DC_*(P^\flat)) \cong \zz$. 

To compute $H_0(\DC_*(P^\flat))$, observe that we may project the kernel of $\epsilon$ onto the subspace generated by any two of the underbridges.
Hence we omit the first row of $\partial_P^\flat$. Then the columns generate an index 3 subgroup of $\zz^2$, showing that $H_0(\DC_*(P^\flat)) \cong \zz/3\zz$.
\end{example}

%%%%%%%%%%%%%%%
\section{Surfaces in the 4-ball and their representations in the 3-sphere}\label{sec:surface}
\begin{figure}[htbp]
\centering
    \begin{subfigure}[t]{0.6\textwidth}
        \centering
        \includegraphics[scale=0.8]{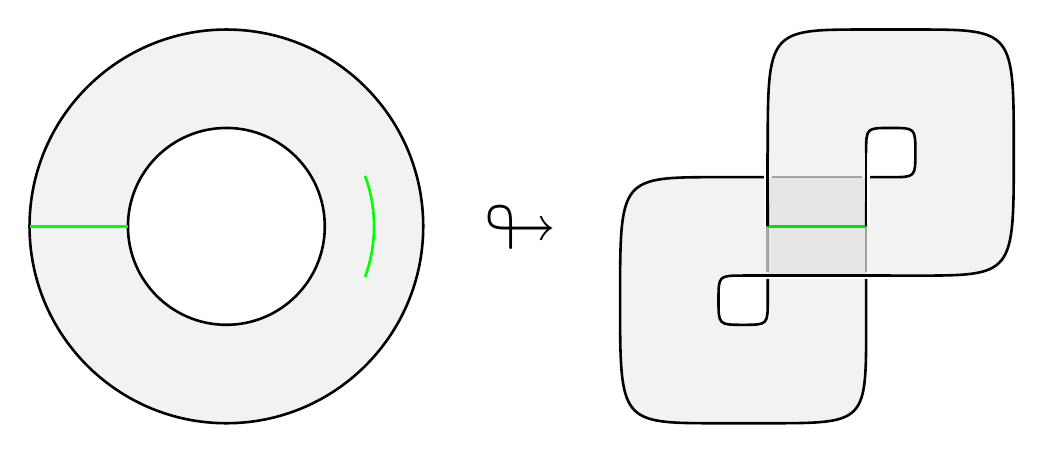}
        \caption{A ribbon-immersed annulus, showing the preimage of the ribbon singularity.}\label{fig:ribbon}
    \end{subfigure}%
    ~ 
    \begin{subfigure}[t]{0.3\textwidth}
        \centering
        \includegraphics[scale=0.8]{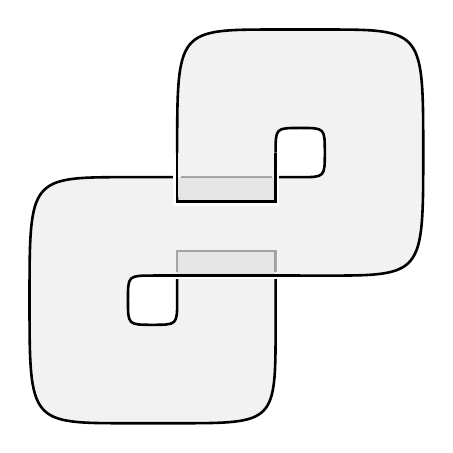}
        \caption{The associated cut surface.}\label{fig:cut}
    \end{subfigure}
\caption
{\bf A ribbon surface $F$ and associated cut surface $F_c$.}
\label{fig:local}
\end{figure}

Let $F \subset B^4$ be a properly embedded compact surface, not necessarily connected or orientable. We will refer to $F$ as a \emph{slice} surface. Denote by $\L\subset S^3$ the link consisting of the boundary components of $F$. We may assume (after an isotopy rel boundary) that the radial distance function $\rho$ in $B^4$ restricts to a Morse function $\rho_F$ on $F$. If $\rho_F$ has no critical points of index 2, then $F$ is called a \emph{ribbon} surface and it admits a ribbon immersion into $S^3$, in which case we denote the image of this immersion by $F_r$. The immersed surface $F_r$ can be described by first choosing pairwise disjoint embeddings of the 0-handles of $F$ into $S^3$ and then connecting them with pairwise disjoint 1-handles that may form ribbon singularities with the images of the 0-handles. Such a surface has a finite number of ribbon singularities as shown in \Cref{fig:ribbon}; the preimage of each consists of two arcs in $F$, one of which is contained in the interior of $F$ (called the \emph{interior arc}) and one which has its endpoints on the boundary of $F$ (called the \emph{properly embedded arc}).  The ribbon-immersed surface $F_r$ is  embedded away from the ribbon singularities where it has two kinds of singular points: interior double points (in the interior of a ribbon singularity) and boundary double points (endpoints of a ribbon singularity). 

Note that $F \subset B^4$ is obtained from $F_r \subset S^3$ by pushing its interior into the interior of $B^4$, where the interior arc of each ribbon singularity is pushed further in than the properly embedded arc. This may be done so that  $\rho_F$ is a Morse function with no maxima. 

Sometimes it will be convenient to convert the immersed surface $F_r$ into an embedded surface by removing a small neighborhood of the properly embedded arc in the preimage of each ribbon singularity; we call this the \emph{cut surface} associated to $F_r$, and denote it by $F_c$; see \Cref{fig:cut}.

For a general slice surface $F \subset B^4$ we may assume that $\rho_F$ is a weakly self-indexing Morse function, i.e., that critical points of higher index have greater radial distance than critical points of lower index.
In particular, we may assume that all minima and saddle points lie in $\rho^{-1}(0,2/3)$, and all maxima in $\rho^{-1}(2/3,1)$. After a further isotopy, supported near the non-critical level $2/3$, we may assume that $F$ is transverse to the sphere of radius $2/3$. Then the sublevel set $\hat{F}=F \cap B^4_{2/3}$ is a properly embedded ribbon surface to which we associate a ribbon-immersed surface $F_r \subset S^3$ as above. We also assume that the radial projection of $F$ to $S^3$ restricts to an embedding on the union of $2$-handles of $F$, that on the interior of each 2-handle this projection is transverse to $F_r$, and that this projection is generic. The last condition restricts possible types of singularities of the projected surface (cf.\ \cite{CKS},  and further detail in \Cref{sec:dbc}).

The boundary $\hat{\L}$ of $\hat{F}$ is the union of two sublinks, $\hat{\L}_0$ and $\hat{\L}_1$.  The first of these is an unlink consisting of those boundary components of $\hat{F}$ that are capped off in $F$ by the 2-handles, and the second corresponds to $\L$, in the sense that a part of the surface $F \smallsetminus \Int \hat{F}$ defines an isotopy between $\hat{\L}_1$ and $\L$.  In particular, the components $L_i$ of $\hat{\L}_0$ bound pairwise disjoint embedded disks $d_i \subset S^3$ (images of 2-handles of $F$) that do not intersect $\hat{\L}_1$, but may intersect the interior of the immersed surface $F_r$.  We call $\hat{\L}_0$ a \emph{separated sublink} of $\hat{\L}$. This yields a 3-dimensional description of the slice surface $F$ as the union, $F_s$, of the ribbon-immersed surface $F_r$ and the disks $d_i$.

%
%The first of these corresponds to $\L$ (in the sense that a part of the surface $F \smallsetminus \Int \hat{F}$ defines an isotopy between $\hat{\L}_0$ and $\L$) and the second  is an unlink consisting of those boundary components of $\hat{F}$ that are capped off in $F$ by the 2-handles. In particular, the components $L_i$ of $\hat{\L}_1$ bound pairwise disjoint embedded disks $d_i \subset S^3$ (images of 2-handles of $F$) that do not intersect $\hat{\L}_0$, but may intersect the interior of the immersed surface $F_r$.  We call $\hat{\L}_1$ a \emph{separated sublink} of $\hat{\L}$. This yields a 3-dimensional description of the slice surface $F$ as the union, $F_s$, of the ribbon-immersed surface $F_r$ and the disks $d_i$.

Note that $F_s$ may not be smoothly embedded along the boundaries of the disks $d_i$. The double points interior to $d_i$ form a 1-manifold that may have closed components and arcs that end on the boundary of $F_r$. These endpoints may either be endpoints of ribbon singularities of the ribbon surface $F_r$ (where two sheets of the projected surface $F_s$ meet transversely) or may indeed be singular points of the projected surface, called \emph{pinch points} or Whitney umbrella singularities which occur when the framing curve from $d_i$ intersects $F_r$; note that the framings determined by $d_i$ and $F_r$ along the common boundary agree modulo 2. The standard model for the Whitney umbrella is given by the solutions of $x^2=y^2z$, with $z\ge0$.  One may then consider the subset with $xy\le 0$ as a part of the ribbon surface and the subset with $xy\ge 0$ as a part of the 2-handle.  The pinch point is at the origin and the double points lie on the $z$-axis.
The last type of singularities in $F_s$ are triple points, where $d_i$ passes through the interior of a ribbon singularity of $F_r$.  

Conversely, a slice surface description $F_s \subset S^3$ determines a slice surface $F \subset B^4$. First its ribbon-immersed subsurface $F_r$ determines a ribbon surface $\hat{F} \subset B^4_{2/3}$. If $\L_0$ is a separated sublink of the boundary of $F_r$ (or equivalently $\hat{F}$), then $\hat{F}$ may be extended to a (possibly closed) slice surface in $B^4$ obtained by capping off the boundary components in $\L_0$.

%%%%%%%%%%%%%%%%%%
\section{Disoriented homology of a ribbon surface}\label{sec:ribbon-disoriented}
The domain of the Gordon-Litherland-type pairing for a ribbon-immersed surface $F_r \subset S^3$ is a subgroup of the first homology group of $F_r$ which we now describe.  We note for emphasis that $F_r$ is the image,  not the domain, of an immersion into the 3-sphere.

\begin{figure}[htbp]
\centering
        \includegraphics{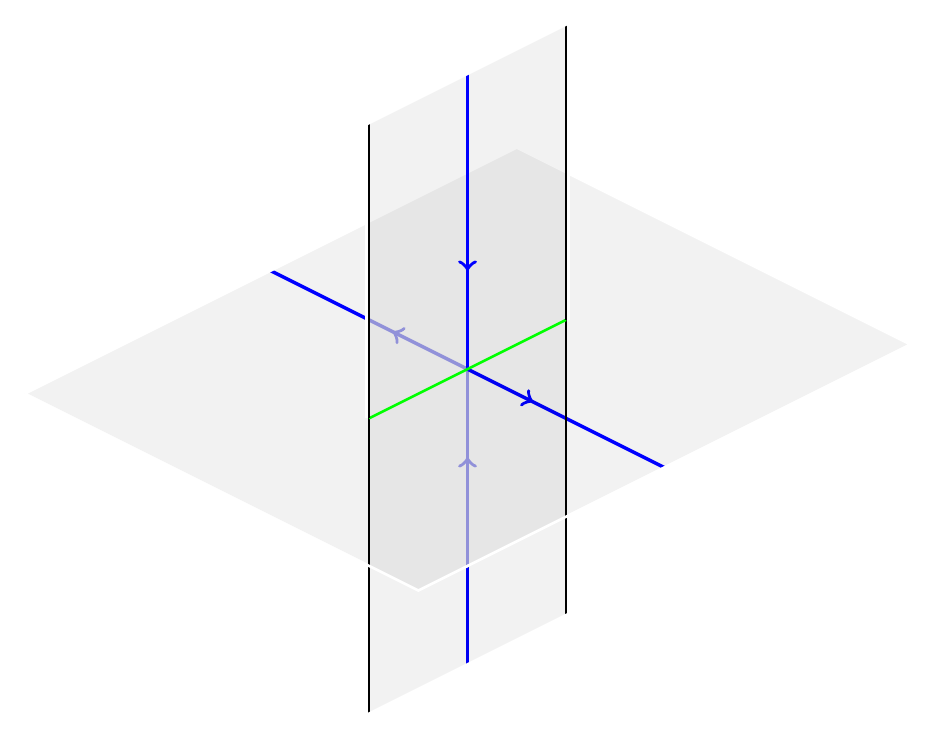} % [scale=0.5]
\caption
{{\bf Local picture near a ribbon singularity.} The local disoriented 1-chain is shown in blue.}
\label{fig:localdis}
\end{figure}

Number the ribbon singularities from $1$ to $k$, and choose coordinates in a cubical ball $B_j$ centred on the $j$th ribbon singularity in $F_r$, with the surface inside the ball consisting of the square $[-2,2]\times[-2,2]\times\{0\}$ in the $(x,y)$-plane and the vertical strip $\{0\}\times[-1,1]\times[-2,2]$ lying in the $(y,z)$-plane.  The \emph{local disoriented $1$-chain} $\ell_j$ associated to the $j$th ribbon singularity is the sum of four oriented line segments in this coordinate patch: two vertical segments, running from $(0,0,\pm2)$ to the origin, and two horizontal segments, running from the origin to $(\pm2,0,0)$. This is sketched in \Cref{fig:localdis}.  
A \emph{disoriented cycle}  is a 1-cycle on $F_r$ of the form
$$a=\sum_{j=1}^k n_j\ell_j + a',$$
where the $n_j$ are integers and $a'$ is a 1-chain supported in the complement of $\Int(B_1\cup\dots\cup B_k)$.
The \emph{(first) disoriented homology group} $DH_1(F_r)$ of the immersed surface $F_r \subset S^3$ is defined to be the subgroup of $H_1(F_r;\zz)$ consisting of classes represented by disoriented cycles.

The disoriented homology group $DH_1(F_r)$ of a ribbon-immersed surface is a free abelian group (as a subgroup of $H_1(F_r)$). In simple situations (such as in the lemma below) it is abstractly isomorphic to the first homology of the underlying surface $F$ which in the presence of ribbon singularities is a smaller group than the homology group of the immersed surface $F_r$. 

\begin{lemma}
\label{lem:DHsimple}
Let $F_r \subset S^3$ be a ribbon-immersed surface with associated ribbon surface $F\subset B^4$.  Suppose that $F$ has a handle decomposition with a single 0-handle and no 2-handles, such that all ribbon singularities are formed by 1-handles passing through the 0-handle; or in other words, the interior arcs are contained in the 0-handle and the properly-embedded arcs are contained in the 1-handles.  Then $DH_1(F_r)$ is isomorphic to $H_1(F;\zz)$.
\end{lemma}

\begin{proof} 
For each 1-handle $h$ of $F$ choose an orientation of its core; these oriented 1-chains, which can be completed to 1-cycles $\gamma_h$ with the addition of oriented arcs in the 0-handle, give a generating set for $H_1(F;\zz)$ as a free abelian group.  We now describe a corresponding set of generators for $DH_1(F_r)$. If a 1-handle $h$ contains no ribbon singularities, let $\alpha_h=\gamma_h$. Otherwise construct a representative for the class $\alpha_h$ by starting with the core of $h$, split into subarcs by the ribbon singularities; orient the first arc arbitrarily and propagate the orientation along the core by changing the orientation of the arc after every ribbon singularity. Let $a_h$ be obtained from this chain by adding appropriately oriented short pairs of arcs in the 0-handle emanating from 
% its intersections with 
the ribbon singularities formed by $h$ as prescribed in \Cref{fig:localdis}. We claim that $a_h$ can be completed to a 1-cycle in $F_r$ by connecting its endpoints with oriented arcs in the 0-handle. Indeed, let $m$ be the number of ribbon singularities along $h$. If $m=2s+1$ is odd, then the endpoints of the core have the same orientation (pointing into or out of the 0-handle) as $2s$ of the other endpoints of oriented arcs comprising $a_h$, whereas the remaining $2(s+1)$ endpoints have the opposite orientation. Similarly, if $m=2s$ is even, the endpoints of the core have opposite orientation, and the other endpoints of oriented subarcs of $a_h$ may be split into two sets of size $2s$ each containing points of one orientation. Hence the endpoints of $a_h$ may be connected up by oriented arcs supported in the 0-handle; we denote the resulting disoriented class by $\alpha_h$. This homology class is well defined, since any two choices of oriented arcs in the 0-handle differ by a trivial cycle.

We now show that any class $\alpha\in DH_1(F_r)$ can be uniquely expressed as a linear combination of $\alpha_h$. By definition, $\alpha$ is represented by a 1-cycle in $F_r$ of the form
$$a=\sum_{j=1}^k n_j\ell_j + a',$$
where $\{\ell_j\}_{j=1}^k$ are the local disoriented 1-chains at ribbon singularities and $a'$ is supported in the complement of the chosen neighborhood of the ribbon singularities. Note that the coefficients $n_j$ are uniquely determined by the homology class $\alpha$ as $F_r$ has the homotopy type of a 1-complex. Consider a 1-handle $h$ of $F$ with $m>0$ ribbon singularities, which we label $j_1,\ldots,j_m$ in the order one encounters them traveling from one end of $h$ to the other. We claim that for all $i<m$, the coefficients satisfy the relation $n_{j_i}+n_{j_{i+1}}=0$. To see this, consider the rectangular part of $h$ between the $i$th and $(i+1)$st ribbon singularity. At the $i$th singularity, the part of $a$ in the rectangle consists of an arc of multiplicity $n_{j_i}$ pointing towards it, and at the $(i+1)$st singularity of an arc of multiplicity $n_{j_{i+1}}$. Since $a$ is a cycle, the sum of the multiplicities of the endpoints of these arcs must be 0 from which the claim follows. We let $n_h=\pm n_{j_i}$ where the sign is chosen so that $n_h\alpha_h$ has the same local multiplicities as $a$ at the ribbon singularities along $h$. Then $\alpha-\sum_h n_h\alpha_h$ is represented by a cycle in the cut surface $F_c$ and hence uniquely expressible as a linear combination of the classes $\alpha_h$ for those 1-handles $h$ that do not contain ribbon singularities.
\end{proof}

\begin{figure}[htbp]
\centering
        \includegraphics[scale=1.2]{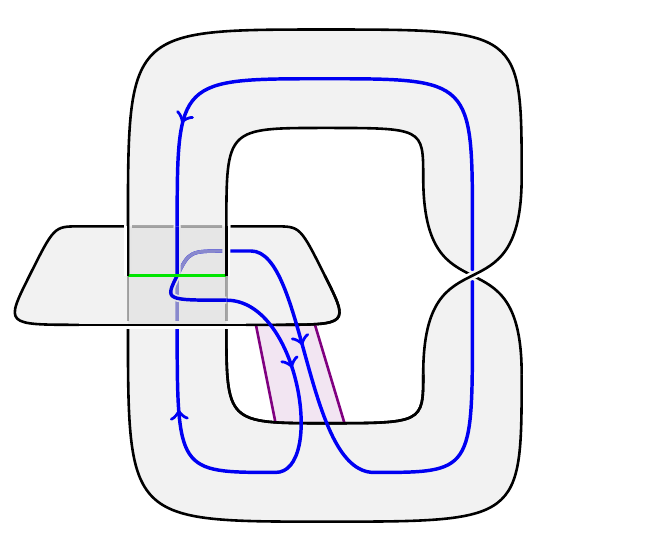}  
\caption
{{\bf A ribbon surface (in grey) with a virtual band (shown in violet).} A generator for the disoriented homology of the virtually banded surface is shown in blue.}
\label{fig:virtualband}
\end{figure}

\Cref{fig:DHgen} shows an example of a ribbon-immersed surface satisfying the hypotheses of the lemma. In general the conclusion of the lemma does not hold as can be seen in the example in Figure \ref{fig:virtualband}.

Suppose that $G_r \subset S^3$ is a ribbon-immersed surface and $F_r$ is a ribbon-immersed subsurface of $G_r$.  If $G_r$ may be obtained from $F_r$ by adding 1-handles (possibly containing ribbon singularities), then the inclusion map of $F_r$ into $G_r$ induces a monomorphism $DH_1(F_r)\to DH_1(G_r)$.

We give another description of $DH_1(F_r)$ which is often easier to work with, and which enables us to also define the 0-dimensional disoriented homology group of $F_r$. 
%{\color{red}We assume that $F_r$ is connected; if not, the same construction may be applied to each of its components.}
If the cut surface $F_c$ is disconnected we attach some embedded 1-handles, which we call \emph{virtual bands}, to $F_r$ to form a new ribbon-immersed surface $G_r \subset S^3$. Denote the collection of virtual bands added to $F_r$ to form $G_r$ by $\calV$. These handles intersect $F_r$ only along their attaching arcs and are pairwise disjoint. They must also satisfy the following conditions:
\begin{enumerate}[(i)]
\item\label{vb1} a virtual band is attached to each component of the cut surface $F_c$ that is not a topological disk with two cuts on the boundary; 
\item\label{vb2} a virtual band is attached to each component of the cut surface $F_c$ containing the interior arc of a ribbon singularity;
\item\label{vb3} the graph $\Pi(\calV)$ with vertices corresponding to the components of $F_c$ to which virtual bands are attached, and edges corresponding to virtual bands is connected. 
\end{enumerate}
For example, virtual bands may be attached to all components of $F_c$ in such a way that the graph in \eqref{vb3} is a tree.

An orientation of a virtual band is an orientation of its core; we fix a choice of orientation for each virtual band. We call $G_r$ a \emph{virtually-banded surface} associated to $F_r$.

To a virtually-banded surface $G_r$ corresponding to $(F_r,\calV)$ we associate a chain complex $\DC_*(F_r,\calV)$ with two nontrivial groups. The group $\DC_1(F_r,\calV)$ is the disoriented homology group $DH_1(G_r)$ of $G_r$; this is typically easier to work with than $DH_1(F_r)$ since it is possible to choose a generating set with at most one generator intersecting each ribbon singularity. The group $\DC_0(F_r,\calV)$ is the free abelian group on $\calV$.  The boundary homomorphism
\begin{align}
\partial_{\calV}:\DC_1(F_r,\calV)&\to \DC_0(F_r,\calV),\label{eq:boundary}\\
[a]&\mapsto\sum_{V \in \calV} \lk(a,K_V)v\notag
\end{align}
is defined as follows. For a virtual band $V\in\calV$ let $K_V$ be the boundary of an oriented disk in $S^3$ whose intersection with $G_r$ is a cocore of $V$, where the orientation of the disk is fixed by the requirement that the intersection number between the disk and the core of $V$ is $+1$. Then the boundary map $\partial_{\calV}$ is given by the linking numbers with the $K_V$'s, or in other words by the signed count of how many times a disoriented homology class passes over each virtual band in the chosen direction.

\begin{proposition}
\label{prop:dishomviaG}  % {lem:dishomviaG}
Let $F_r \subset S^3$ be a ribbon-immersed surface and let $G_r$ be a virtually-banded surface corresponding to $(F_r,\calV)$ as above. Then $H_*(\DC_*(F_r,\calV))$ is (up to isomorphism) independent of the choices in the construction of $G_r$, and the inclusion of $F_r$ in $G_r$ induces a canonical isomorphism
$$DH_1(F_r)\cong H_1(\DC_*(F_r,\calV)).$$
We call the homology of the chain complex $\DC_*(F_r,\calV)$ the \emph{disoriented homology} of $F_r$, denoted by $DH_*(F_r)$.
\end{proposition}

The relation of $DH_0(F_r)$ to the 4-dimensional description will be made apparent in \Cref{sec:dbc4}.

\begin{proof}  
%{\color{red} We assume $F_r$ and consequently $G_r$ is connected.}
We construct a handle decomposition of $G_r$ without 2-handles and with a single 0-handle containing all the interior arcs of ribbon singularities as follows. Start with a handle decomposition without 2-handles of (the underlying surface of) $F_r$, so that each component of the cut surface $F_c$ to which a virtual band is attached contains a single 0-handle, and there are no other 0-handles. We may assume that the virtual bands are attached to the 0-handles, that the interior arcs of the ribbon singularities are contained in the interiors of the 0-handles, and that the properly embedded arcs of ribbon singularities are contained in the 1-handles. 

Recall the graph $\Pi(\calV)$ in condition \eqref{vb3} governing the attachment of virtual bands. Since this graph is connected, we may choose $\calV_0 \subset \calV$ so that  the graph $\Pi(\calV_0)$ is a maximal tree of $\Pi(\calV)$.  Then the union of the 0-handles in $F_r$ and the virtual bands in $\calV_0$ forms a single 0-handle in the decomposition of $G_r$ that contains all the interior arcs. Hence by Lemma \ref{lem:DHsimple} the group $\DC_1(F_r,\calV)=DH_1(G_r)$ is isomorphic to the free abelian group with one generator for each 1-handle of $G_r$; these are 1-handles of $F_r$ and virtual bands not in $\calV_0$.

Since the ribbon-immersed surface $G_r$ is obtained from $F_r$ by adding embedded 1-handles, the inclusion of $F_r$ into $G_r$ induces an inclusion of $DH_1(F_r)$ as a subgroup of $DH_1(G_r)$. Note that disoriented 1-cycles in $F_r$ do not intersect the virtual bands, so these are in the kernel of $\partial_{\calV}$. On the other hand, if a disoriented 1-cycle in $G_r$ is in the kernel of $\partial_{\calV}$, then it is homologous to a disoriented 1-cycle in $F_r$ by a homology supported in the virtual bands of $G_r$. This proves that the inclusion $F_r \hookrightarrow G_r$ induces a canonical isomorphism $H_1(\DC_*(F_r,\calV))\cong DH_1(F_r)$.

To prove independence of $H_0(\DC_*(F_r,\calV))$ from the choices made in the construction first note that adding a virtual band $V$ to $\calV$ and hence to $G_r$ subject to the above conditions yields a new surface $G_r'$ corresponding to $\calV'=\calV \cup\{V\}$ with isomorphic $H_0$. Indeed, the chain complex $\DC_*(F_r,\calV')$ is obtained from $\DC_*(F_r,\calV)$ by adding a generator to each of its groups: 
$$ \DC_0(F_r,\calV')=\DC_0(F_r,\calV)\oplus \zz v, \ \ \DC_1(F_r,\calV')=\DC_1(F_r,\calV)\oplus  \zz \alpha, $$
where $\alpha$ is represented by a 1-cycle in $G_r'$ that passes over the virtual band $V$ geometrically once. This shows that the inclusion of $\DC_*(F_r,\calV)$ into $\DC_*(F_r,\calV')$ is a chain equivalence. Note that by condition \eqref{vb3} above at least one of the two (possibly the same) components of $F_c$ that $V$ connects is in $G_r$ already connected to a virtual band, and hence to the 0-handle of $G_r$. If this holds for both components, then we take $\calV_0'=\calV_0$ and $V$ becomes a 1-handle of $G_r'$.  The core of $V$ may be completed in the 0-handle of $G_r$ to a 1-cycle in $G_r'$ proving the claim in this case. Otherwise $V$ is connected to a component $A$ that in $G_r$ does not have a virtual band attached to it, thus $A$ is a part of a 1-handle $h$ in $G_r$ between two ribbon singularities. We take $\calV_0'=\calV_0 \cup \{V\}$ and change the handle decomposition of $F_r$ by introducing another 0-handle in $A$; this splits $h$ into two 1-handles. Recall from the proof of \Cref{lem:DHsimple} that to $h$ corresponds a generator $\alpha_h$ in $DH_1(G_r)$ constructed from the chain $a_h$. Then one half of $a_h$ (corresponding to one of the new 1-handles) along with the core of $V$ can be as in the proof of that lemma completed to a cycle in $G_r'$ defining the class $\alpha$. Since the coefficient of $v$ in $\partial_{\calV'}(\alpha)$ is $\pm1$, the claim follows.
%boundary homomorphism $\partial_{\calV'}$ maps $\alpha$ to $\pm v$

It follows from the previous paragraph that we may assume virtual bands in $\calV$ are attached to all components of the cut surface $F_c$ and that the corresponding graph $\Pi(\calV)$ is a tree. We now verify that $H_0$ agrees for such choices of collections of virtual bands. Let $\calV_1$ and $\calV_2$ be two such collections and denote by $G_r^1$ and $G_r^2$ the corresponding virtually banded surfaces. Then $G_r^1$ can be transformed into $G_r^2$ by a sequence of steps where in each step a virtual band $V_1 \in \calV_1$ is replaced by a virtual band $V_2 \in \calV_2$;  we may assume by an isotopy that $V_2$ is disjoint from virtual bands in $\calV_1$. By induction we assume there is just one such step, so that $\calV_1 \smallsetminus \{V_1\} = \calV_2 \smallsetminus \{V_2\}$.  Adding $V_2$ to $G_r^1$ gives a new larger surface $\widehat G$ whose graph $\Pi(\calV_1 \cup \{V_2\})$ contains a cycle that includes $V_1$ and $V_2$. This graph cycle gives rise to a cycle in the first homology group of $\widehat G$.  The homology class of this cycle (oriented consistently with $V_1$) is represented by a 1-chain $b$, which may be assumed to be 
disjoint from all ribbon singularities. Then for any class $\alpha=[a] \in \DC_1(F_r, \calV_1)$ we let $\varphi(\alpha)=[a - \lk(a,K_{V_1})b]$; clearly the cycle on the right may be represented in $G_r^2$. On $\DC_0(F_r, \calV_1)$ we let $\varphi$ act as the identity except that it sends $v_1$ to $-\sum_{V\in \calV_2} \lk(b,K_V)v$. Clearly $\varphi$ is a well defined isomorphism of the chain complexes.
\end{proof}

\begin{example}
Consider the ribbon-immersed surface $F_r \subset S^3$ shown in \Cref{fig:virtualband}. The disoriented chain complex for the indicated virtual band is 
$$ \DC_1 \cong \zz,\quad \DC_0 \cong \zz,$$
where the boundary homomorphism is multiplication by $2$. Hence the disoriented homology is  
$$DH_1(F_r)=0,\quad DH_0(F_r) \cong \zz/2\zz.$$
Let $F \subset B^4$ be a properly embedded surface obtained by pushing the interior of $F_r$ into the 4-ball and let $X$ be the branched double cover of the 4-ball with branch set $F$. Then according to \Cref{thm:4d-homology} the reduced homology of $X$ is nontrivial only in dimension $1$ and
$$H_1(X;\zz) \cong \zz/2\zz.$$
\end{example}

%%%%%%%%%%%%%%%%%%%%%%%
\section{Disoriented homology of a slice surface}\label{sec:slice-disoriented}
A slice surface $F \subset B^4$ can be described (as in Section \ref{sec:surface}) by $F_s \subset S^3$ which consists of a ribbon-immersed surface $F_r \subset  S^3$ along with a separated sublink $\L_0=\{L_1,\ldots,L_m\}$ of its boundary. Boundary components in $\L_0$ bound pairwise disjoint disks $d_i \subset S^3$ that do not intersect the rest of the boundary. Choose disjoint small closed regular neighborhoods $N_i$ of $d_i$. We assume $N_i$ is small enough so that $N_i \cap F_r$ is a regular neighborhood of $d_i \cap F_r$ and that the boundary spheres $S_i$ of $N_i$ intersect $F_r$ transversely in its interior. Then the intersection $S_i \cap F_r$ is a 4-valent graph whose vertices are the intersections of $S_i$ with the ribbon singularities. 

\begin{lemma}
With the notation as above, the intersection $S_i \cap F_r$ determines a disoriented 1-cycle $b_i$ and hence a homology class $\beta_i$ in $DH_1(F_r)$, well defined up to sign. 
\end{lemma}

\begin{proof}
Color the faces\footnote{The components of the complement.} of the graph $S_i \cap F_r$ on $S_i$ in a chessboard fashion. A choice of orientation of the sphere induces orientations of the faces. Orient the edges of the graph consistently with the black faces, as  shown in Figure \ref{fig:discyc}. Then the orientation on the graph is consistent with it representing a disoriented 1-cycle $b_i$ in $F_r$.

\begin{figure}[htbp]
\centering
\includegraphics[scale=1]{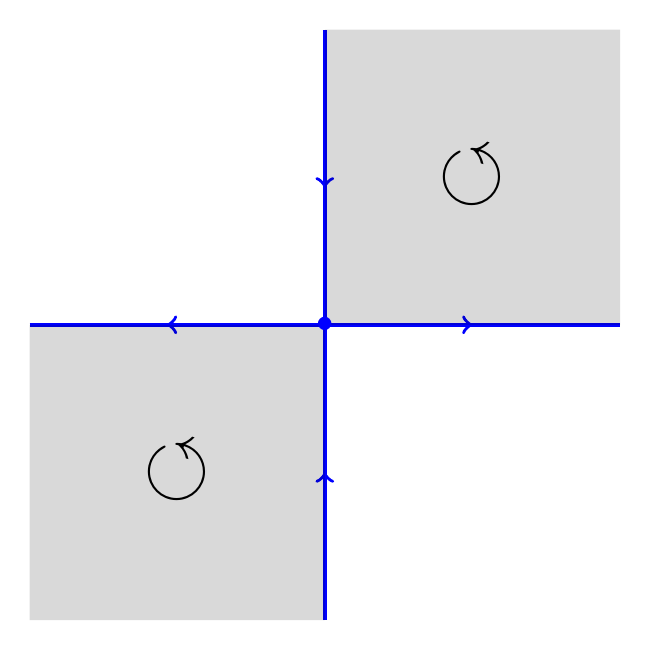}
\caption
{{\bf The disoriented cycle $b_i$.} Orienting the edges of the graph $S_i \cap F_r$ as the boundary of the black faces yields a disoriented homology class. Recall that vertices of the graph come from ribbon singularities.}
\label{fig:discyc}
\end{figure}

If $N_i'$ is another small neighborhood of $D_i$ as above, then the corresponding disoriented cycle $b_i'$ is homologous to $\pm b_i$ since there is a homotopy transforming one into the other.
\end{proof}

Note that in fact $L_i$ uniquely determines the class $\beta_i$. The sphere $S_i$ can be chosen as any separating sphere for this component of $\partial F_r$. Any two such spheres are isotopic in the complement of $\partial F_r$ and thus their intersections with $F_r$ determine the same disoriented homology class (up to sign).

Let $G_r$ be any virtually-banded surface associated to $F_r$ through a choice of virtual bands $\calV$. Together with the link $\L_0$ it determines a disoriented chain complex $\DC_*(F_r,\calV,\L_0)$ as follows. We let $\DC_k(F_r,\calV,\L_0)=\DC_k(F_r,\calV)$ for $k=0,1$, and extend this complex to include another group, $\DC_2(F_r,\calV,\L_0)$, which is the free abelian group with basis the disks $d_i$. The boundary homomorphism $\partial_{\calV} \colon \DC_2(F_r,\calV,\L_0) \to \DC_1(F_r,\calV)$ sends $d_i$ to $\beta_i$. Since the support of $\beta_i$ lies in $F_r$, it follows that $(\DC_*(F_r,\calV,\L_0),\partial_{\calV})$ is indeed a chain complex. 

\begin{definition}\label{def:DHv1}
We call the homology of the complex $(\DC_*(F_r,\calV,\L_0),\partial_{\calV})$ the \emph{disoriented homology} of the slice surface description $F_s$ and denote it by $DH_*(F_s)$. 
\end{definition}

It is clear from the case of ribbon surfaces that the resulting homology is independent of the choice of virtual bands.

%%% Cellular disoriented homology

We give yet another description of the disoriented homology of a slice surface that is defined in terms of a handle decomposition of the surface. This is analogous to the disoriented homology of a link and provides a convenient way of identification with the homology of the double branched cover of the 4-ball.

A handle decomposition of $F$ determines a ribbon subsurface of $F$. We will refer to the images of the handles of $F$ in the corresponding ribbon-immersed surface $F_r$ also as handles. We assume that all the ribbon singularities are formed by 1-handles passing through the 0-handles of $F_r$. The handle decomposition also determines a separated sublink $\L_0$ of the boundary of $F_r$ whose components bound disks $d_i$ (the 2-handles). Given this, choose for each 1-handle $h_j$ of $F_r$ a \emph{disorientation} of its core, i.e., orient the arcs into which ribbon singularities split the core in such a way that any two consecutive arcs have opposite orientations. Denote the disoriented core of $h_j$ by $c_j$. Let $\Gamma_i$ be the intersection of the disk $d_i$ with $F_r$; we assume that this intersection is transverse in the interior of $d_i$. Then $\Gamma_i$ is a graph that contains all of $\partial d_i$ and whose interior vertices are 4-valent corresponding to ribbon singularities. Its vertices on the boundary are 3-valent and correspond to pinch points or ribbon singularities. Choose a chessboard coloring of the faces of $\Gamma_i$ on $d_i$. Then orienting all the black faces consistently with one orientation of the disk $d_i$ and giving all the white faces the opposite orientation determines a \emph{disorientation} of the 2-handle $d_i$ -- we denote the disk along with the chosen disorientation by $d_i^\flat$. We denote this slice surface description of $F$ with a chosen handle decomposition of $F_r$ and chosen disorientations of its 1- and 2-handles as described above by $F_s^\flat$. 

The disoriented chain complex for $F_s^\flat$ is given as follows:
\begin{itemize}
\item $\DC_0(F_s^\flat)$ is the free abelian group generated by the 0-handles;
\item $\DC_1(F_s^\flat)$ is the free abelian group generated by the disoriented cores of the 1-handles;
\item $\DC_2(F_s^\flat)$ is the free abelian group generated by the disoriented 2-handles.
\end{itemize}
The boundary homomorphism $\partial_1^\flat : \DC_1(F_s^\flat) \to \DC_0(F_s^\flat)$ is given by the signed count of the number of times a given disoriented core points into (positive contribution) or away from (negative contribution) a 0-handle; note that the contribution at each ribbon singularity is $\pm 2$ times the zero handle containing the interior arc of the singularity. To define the boundary homomorphism $\partial_2^\flat : \DC_2(F_s^\flat) \to \DC_1(F_s^\flat)$, orient the edges of $\Gamma_i$ as the boundary of the black regions in $d_i^\flat$. This data determines a disoriented homology class $\beta_i^\flat=[b_i^\flat]$ in the ribbon-immersed surface $F_r$ by letting $b_i^\flat$ be the sum of the boundaries of the oriented faces of $\Gamma_i$. Hence $b_i^\flat$ is the linear combination of the oriented edges of $\Gamma_i$, where the edges lying in $\partial d_i$ have multiplicity 1 and the interior edges have multiplicity 2. For each 1-handle $h_j$ of $F_r$ we count how many times $b_i^\flat$ passes over it as follows: choose an orientation of $h_j$
and orient one of its attaching arcs $a_j$ so that the intersection number of $a_j$ and $c_j$ equals 1, $a_j \cdot c_j=1$. Then the coefficient of $c_j$ in $\partial_2^\flat(d_i^\flat)$  is equal to the intersection number $a_j \cdot b_i^\flat$. 

\begin{definition}\label{def:DHv2}
We call the augmented chain complex $(\DC_*(F_s^\flat),\partial_*^\flat)$, where the augmentation homomorphism $\partial_0^\flat=\varepsilon : \DC_0(F_s^\flat) \to \DC_{-1}(F_s^\flat)=\zz$ sends each 0-handle to 1, the \emph{cellular disoriented complex} of the slice surface description $F_s^\flat$. 
\end{definition}

\begin{proposition}\label{prop:DHiso}
The homology of the cellular disoriented complex $(\DC_*(F_s^\flat),\partial_*^\flat)$ is isomorphic to the disoriented homology of $F_s$. 
\end{proposition}

\begin{proof}
We will construct a chain equivalence $f_* : \DC_*(F_s^\flat) \to \DC_*(F_r,\calV,\L_0)$ for a particular choice of a virtually banded surface $G_r$, determined by a collection of virtual bands $\calV$ for $F_r$. Choose a 0-handle $m_0$ of $F_r$ and connect this 0-handle to every other 0-handle $m_i$, $i \ge 1$, by a virtual band $V_i$. Orient virtual bands so that they point to $m_0$. Let $f_0$ be given by
$$f_0(m_i)=v_i,\ i \ge 1, \quad f_0(m_0)=0.$$

By \Cref{lem:DHsimple}, $\DC_1(F_r,\calV,\L_0)=DH_1(G_r)$ is generated by elements corresponding to 1-handles of $F_r$. In fact, a generator $\alpha_j$ corresponding to a 1-handle $h_j$ may be obtained by completing the disoriented core $c_j$ of $h_j$ to a disoriented 1-cycle $\bar c_j$ in $G_r$. This defines the homomorphism $f_1$:
$$f_1(c_j)=\alpha_j=[\bar c_j].$$
Finally, $f_2$ is given by sending each disoriented 2-handle $d_i^\flat$ to the disk $d_i$.
$$
\begin{CD}
\DC_2(F_s^\flat) @>{\partial_2^\flat}>> \DC_1(F_s^\flat)  @>{\partial_1^\flat}>> \DC_0(F_s^\flat) @>{\varepsilon}>> \zz   \\
@V{\cong}V{f_2}V @V{\cong}V{f_1}V @VV{f_0}V  @VV{f_{-1}}V  \\
\DC_2(F_r,\calV,\L_0) @>{\partial_{\calV}}>>  \DC_1(F_r,\calV,\L_0) @>{\partial_{\calV}}>{\phantom{\partial}}>  \DC_0(F_r,\calV,\L_0) @>{\partial_{\calV}}>>  0 
\end{CD}
$$

We verify that $f_*$ is a chain map. Note that the (algebraic count of the) number of times $\bar c_j$ goes over the virtual band $V_i$ connecting $m_i$ to $m_0$ is the same as the coefficient of $m_i$ in $\partial_1^\flat c_j$, hence
$$\partial_{\calV} \circ f_1=f_0 \circ \partial_1^\flat.$$
To show the commutativity of the left square, we need to see that for each disk $d_i$, the resulting 1-cycles $b_i$ and $f_1(b_i^\flat)$ give the same element of disoriented homology (up to sign).
The sphere $S_i$ is a double push-off of the disk $d_i$ and therefore a chessboard coloring of $d_i$ determines a chessboard coloring of $S_i$ by changing all the colors on one of the hemispheres. This orients the two edges of $S_i \cap F_r$ corresponding to a given interior edge of $\Gamma_i$ consistently. A homotopy collapsing the sphere $S_i$ onto the disk $d_i$ now induces a homology between $b_i$ and $\pm b_i^\flat$, as elements of the first homology group of the immersed surface. We choose the sign of $b_i$ so that 
$$\partial_{\calV} \circ f_2=f_1 \circ \partial_2^\flat$$
holds.

We claim that $f_*$ has a chain homotopy inverse $g_* : \DC_*(F_r,\calV,\L_0) \to \DC_*(F_s^\flat)$, where $g_i=f_i^{-1}$ for $i=1,2$, and $g_0$ is given by $g_0(v_i)=m_i-m_0$. Clearly $g_*$ is also a chain map: the commutativity of the left and right squares is clear, for the middle square it follows from the argument for $f_*$ and the choice of $g_0$.
Note that also $f_0 \circ g_0=id$, whereas $g_0 ( f_0(m_i))=m_i-m_0$ for $i \ge 1$, and $g_0 ( f_0(m_0))=0$. Hence a chain homotopy between $id$ and $g_* \circ f_*$ is given by $H : \DC_*(F_s^\flat) \to \DC_{*+1}(F_s^\flat)$ whose only nontrivial component is $H_{-1}$ which sends $1 \in \zz$ to $m_0$:
\begin{equation*}
id-g_0 \circ f_0= H_{-1} \circ \varepsilon. \qedhere
\end{equation*}
\end{proof}

%%%%%%%%%%%%%%%
\section{The Gordon-Litherland type pairing on the disoriented homology group}
\label{sec:pairing}
Let $F_r \subset  S^3$ be a ribbon-immersed surface. Given two disoriented homology classes $\alpha$, $\beta \in DH_1(F_r)$ represented by disoriented cycles $a$ and $b$ we wish to follow Gordon and Litherland \cite{gl}, and define the pairing of $\alpha$ and $\beta$ to be the linking number of $a$ and $\tau b$,
where $\tau b$ is obtained by pushing $b$ off $F_r$ in the normal direction on both sides.  

\begin{figure}[htbp]
\centering
        \includegraphics{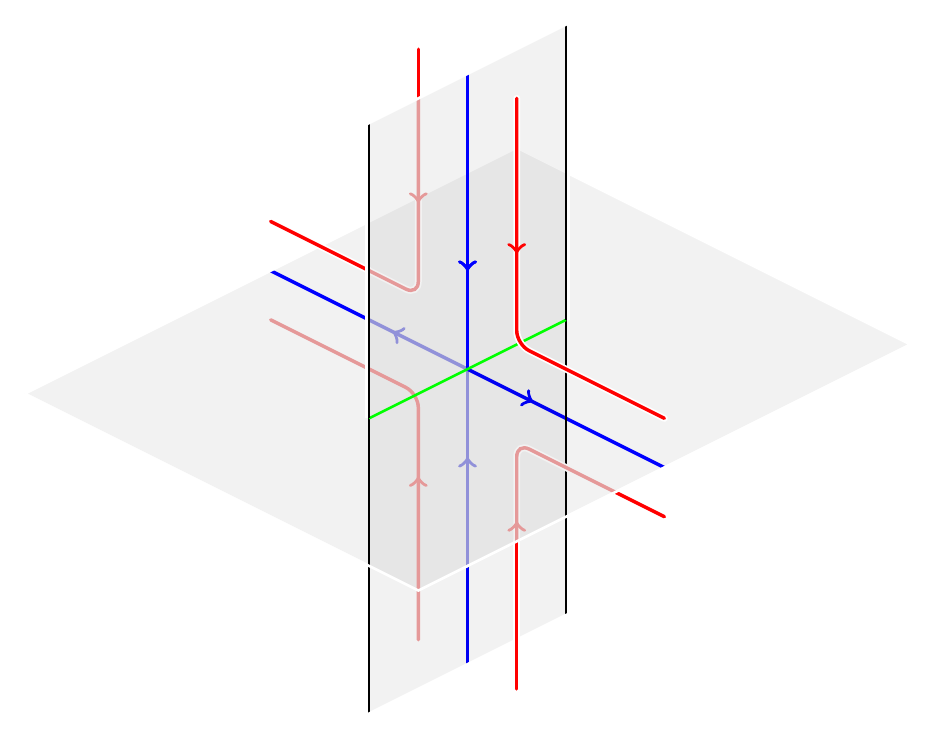} % [scale=0.5]
\caption
{{\bf The double push-off near a ribbon singularity.} The local disoriented 1-chain is shown in blue, with its double push-off in red.}
\label{fig:pushoff}
\end{figure}

Of course, we need to take care in defining $\tau b$ in the vicinity of a ribbon singularity.  
Recall that $b$ is represented by a 1-chain on $F_r$ whose support near each ribbon singularity is an integer multiple of the local disoriented 1-chain $\ell_j$ shown in \Cref{fig:localdis}.  We take coordinates in a  ball neighbourhood $B_j$ of each ribbon singularity as before.  The push-off $\tau \ell_j$ of $\ell_j$ then consists of two disjoint oriented line segments in each of the planes $z=\pm1$ and $x=\pm1$.  The starting points of the segments in the plane $z=1$ are $(\pm1,0,1)$ and the endpoints are $(\pm2,0,1)$.  We then take the vertical translates of these two segments in the plane $z=-1$.  Similarly  the segments in the plane $x=1$ go from $(1,0,\pm1)$ to $(1,0,\pm2)$, and we take horizontal translates of these in the plane $x=-1$.
Away from the ribbon singularity we take normal pushoffs on either side of $F_r$ as usual, chosen to match up with (the given multiple of) $\tau \ell_j$.  The result is a (singular if $\max |n_j|>1$) oriented link $\tau b$ in $S^3\smallsetminus F_r$, as illustrated in \Cref{fig:pushoff}.

The Gordon-Litherland type form for the ribbon-immersed surface $F_r \subset S^3$ is now defined to be
$$\lambda_{F_r}(\alpha,\beta)=\lk(a,\tau b).$$
Note that $a$ and $b$ in the above formula may be singular; see \Cref{subsec:lk} for discussion of linking numbers in this case.

\begin{example}
\label{eg:DHgen}
One may check that the square $\lambda_{F_r}(\alpha,\alpha)$ of the generator shown in \Cref{fig:DHgen} is $6$, which agrees with the determinant of the boundary of the given surface.  
\end{example}

\begin{proposition}
For a ribbon-immersed surface $F_r$, $\lambda_{F_r}$ is a well-defined symmetric bilinear form on $DH_1(F_r)$. Moreover, if $G_r$ is a ribbon-immersed surface obtained from $F_r$ by adding 1-handles, then the restriction of $\lambda_{G_r}$ to the disoriented homology of $F_r$ agrees with $\lambda_{F_r}$.
\end{proposition}

\begin{proof}
Since the linking number of two disjoint cycles $a$ and $\tau b$ depends only on the homology classes of the cycles, and since a homology between $b$ and $b'$ in $F_r$ naturally gives rise to a homology between $\tau b$ and $\tau b'$ in the complement of $F_r$, it follows that $\lambda_{F_r}$ is well-defined. That $\lambda_{F_r}$ is symmetric follows similarly as in the case of embedded surfaces \cite{gl}. Let $N$ be the immersed normal $B^1$-bundle of $F_r$ in $S^3$. Self-intersections of $N$ are cubes located at ribbon singularities of $F_r$. Denote by $\partial ' N$ the part of the boundary of $N$ that comes from the $S^0$-bundle; we smooth the corners in  $\partial ' N$ along the edges of the cubes at ribbon singularities. Let the positive normal direction to $\partial ' N$ be given by the outward pointing normal and for any 1-cycle $c$ on $\partial ' N$ denote by $c^+$ a nearby pushoff of $c$ in this direction and by $c^-$ a nearby pushoff of $c$ in the opposite direction. Note that $\tau a$ may be viewed as a 1-cycle in $\partial ' N$ and that it is homologous to $2a$ in $N$. Then
$$\lk(a,\tau b)= \lk(a,\tau b^+)=\lk(\tau a,\tau b^+)/2$$
and therefore
\begin{align*}
2\big(\lambda_{F_r}([a],[b])-\lambda_{F_r}([b],[a])\big) &= \lk(\tau a,\tau b^+)- \lk(\tau b,\tau a^+)\\
&= \lk(\tau a, \tau b^+ - \tau b^-) = \tau a \cdot B,
\end{align*}
where $B$ is the 2-chain with boundary $\tau b^+ - \tau b^-$, obtained by restricting the normal $B^1$-bundle of $\partial ' N$ to $\tau b$. 
Note that each intersection point $x$ between $a$ and $b$ in $F_r$ gives rise to a pair of intersection points $\tau x$ between $\tau a$ and $\tau b$ in $\partial ' N$ at which the orientations of the normal to $\partial ' N$ are opposite. In other words, the two patches of $\partial ' N$ at the points $\tau x$ have opposite orientations. Hence the local intersection numbers at the two points in $\tau x$ are of the opposite sign and the intersection number above vanishes.

The last claim of the proposition is clear from the definition of the pairing.
\end{proof}

Consider now a description $F_s \subset S^3$ of a slice surface, consisting of a ribbon-immersed surface $F_r \subset  S^3$ and a separated sublink $\L_0$ of its boundary. The following lemma shows that the form $\lambda_{F_r}$ induces a well-defined symmetric bilinear form $\lambda_{F_s}$ on $DH_1(F_s)$. Recall that to any component $L_i$ of $\L_0$ we associate a class $\beta_i \in DH_1(F_r)$ represented by a disoriented 1-chain $b_i$ whose support is the intersection of a separating sphere $S_i$ for $L_i$ with $F_r$.

\begin{lemma}
With the notation as above, $\lambda_{F_r}(\alpha,\beta_i)=0$ for any $\alpha \in DH_1(F_r)$.
\end{lemma}

\begin{proof}
Since the sphere $S_i$ is transverse to $F_r$ we may take the double pushoff $\tau b_i$ to be the boundary of a bicollar neighborhood of $b_i$ in $S_i$. Recall that $b_i$ is oriented consistently with the black regions in a chessboard coloring of its complement. The complement of the open bi-collar is a union of disks and if we let $c$ be the 2-chain which is given by the sum of all the black disks minus the sum of all the white disks, then $\tau b_i$ is the boundary of this 2-chain. Since $c$ does not intersect $F_r$, it follows that the linking number of $\tau b_i$ with any (disoriented) 1-cycle on $F_r$ is zero.
\end{proof}

\begin{example}[The positive unknotted real projective plane]\label{ex:proj-homology}
We compute the disoriented homology and the GL-pairing of the unknotted real projective plane $P=\rr\pp^2$ in $B^4$ with radial projection $P_s$ given on the left of \Cref{fig:proj-homology}.
\begin{figure}[htbp]
\centering
\includegraphics[scale=1.2]{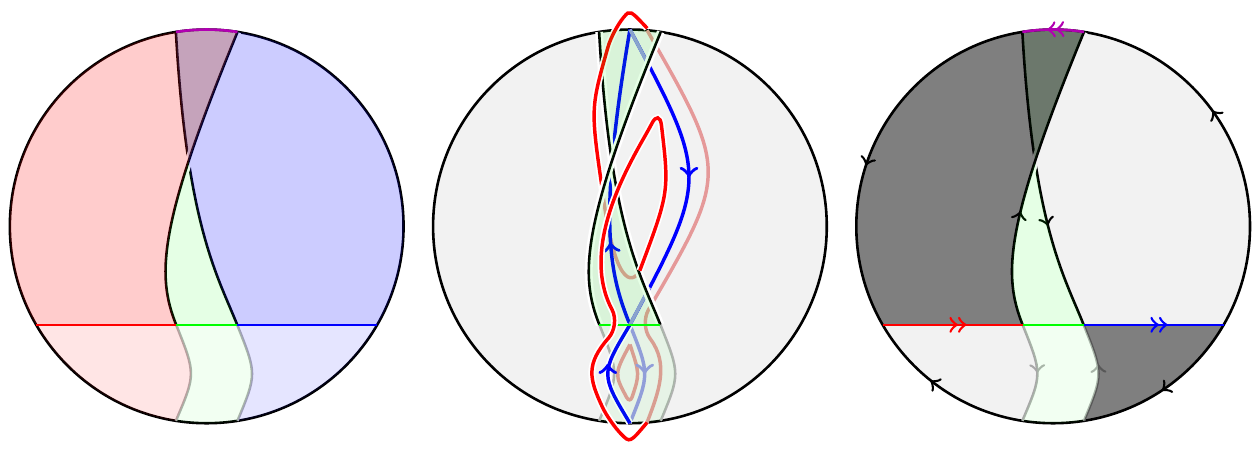}  
\caption
{{\bf The real projective plane.} The left picture shows the radial projection $P_s$ of $P$: the round disk is the 0-handle, the green band the 1-handle, and the red and blue disks combine to give the 2-handle. The middle picture shows a generator for the first disoriented homology (in blue) and its pushoff (in red). The right picture gives a disorientation of the 2-handle and the resulting cycle $b^\flat$.}
\label{fig:proj-homology}
\end{figure}
Denote the 0-handle of $P$ by $m$, its 1-handle by $h$ and 2-handle by $d$. Let $c_h$ be the disoriented core of $h$, given by the part of the blue generator in the middle picture on \Cref{fig:proj-homology} lying on $h$. Choose a disorientation $d^\flat$ of $d$ as in the right picture on \Cref{fig:proj-homology}.  Then the cellular disoriented chains of $P$ are
$$\DC_0(P_s^\flat)=\zz m,\quad \DC_1(P_s^\flat)=\zz c_h,\quad \DC_2(P_s^\flat)=\zz d^\flat.$$
The boundary homomorphism on $\DC_1$ is trivial as there is only one 0-handle. Also the boundary homomorphism on $\DC_2$ is trivial as can be seen from the right picture on \Cref{fig:proj-homology} since the two arcs of the boundary cycle $b^\flat$ of $d^\flat$ have opposite disorientations. Taking into account the augmentation homomorphism it follows that 
$$DH_0(P_s^\flat)=0,\quad DH_1(P_s^\flat)\cong \zz,\quad DH_2(P_s^\flat)\cong\zz.$$
The self-pairing of the generator of  $DH_1(P_s^\flat)$ is equal to $+1$ as can be seen from the middle picture on \Cref{fig:proj-homology}, since the linking number between the disoriented cycle in blue and its pushoff in red is equal to $+1$.

Note that if one changes the crossing in the projection of $P$ one obtains the negative unknotted projective plane; one sees immediately that the disoriented homology groups do not change but the sign of the pairing on $DH_1$ switches to negative.

We return to this example in \Cref{ex:proj-Kirby} where we exhibit a Kirby diagram of the branched double cover of the 4-ball with branch set $P$.
\end{example}

%%%
\subsection{Remarks on linking numbers}
\label{subsec:lk}
Given two disjoint oriented knots $K$ and $K'$ in $\rr^3$, represented as smooth maps from $S^1$ to $\rr^3$,  their linking number $\lk(K,K')$ may be defined as the degree of the map
$$(u,v)\mapsto\frac{K(u)-K'(v)}{|K(u)-K'(v)|}$$
from $S^1\times S^1$ to $S^2$. This implies that the linking number is an invariant of homotopy classes of maps with disjoint images and also that it is symmetric. The linking number is then extended to links by requiring it to be bilinear: if $L=K_1\cup\dots\cup K_m$ and $L'=K'_1\cup\dots\cup K'_n$ are disjoint oriented links, then 
$$\lk(L,L')=\sum_{i,j}\lk(K_i,K'_j).$$

Alternatively, $\lk(K,K')$ may be defined as the multiple of the homology class determined by $K'$ in $H_1(\rr^3 \smallsetminus K;\zz) \cong \zz$, where the generator is a positively oriented meridian of $K$. So in fact the linking number depends only on the homology class of $K'$ in the complement of $K$. To explicitly compute $\lk(K,K')$, one usually relies on combinatorial description via diagrams: starting with a diagram of $K \cup K'$, assign to each crossing $c$ between $K$ and $K'$ a sign $\veps_c \in \{\pm 1\}$, where $\veps_c=1$ if a bug travelling along the overcrossing arc in the chosen direction sees the undercrossing arc oriented from right to left.  Then
$$\lk(K,K') = \frac12 \sum_c \veps_c;$$
if one counts only overcrossings of one knot over the other, the same formula without the half applies.

\begin{figure}[htbp]
\centering
\includegraphics[scale=1.2]{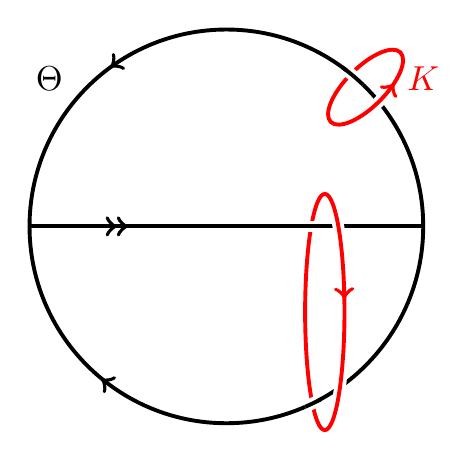}  
\caption
{\bf The linking number of the embedded $\Theta$-graph with the knot $K$ is $+1$.  Two isotopy representatives of $K$ are shown.}
\label{fig:thetalinking}
\end{figure}

As pointed out in \cite{rolfsen}, the definition allows for each of $L$ and $L'$ to be singular, as long as they are disjoint. In fact, $L$ and $L'$ may be any two disjoint 1-cycles. The case of interest to us is when $L$ and/or $L'$ is an embedded graph with oriented edges, with nonnegative integer multiplicities associated to each edge, in such a way that the signed weighted count of edges at each vertex (inward minus outward) is zero.   An example is shown in \Cref{fig:thetalinking}; any such graph has an interpretation as a singular link in which the multiplicity of an edge is the signed number of times it is traversed by the components of the link. Clearly one may apply the above combinatorial formula to compute the linking number of such objects.

%%%%%%%%%%%%%%%%%%%%%%%%%%%%%%%%%%%%%%%%%%%%%%%%%%%%%%%%%%%%%%%%%%%%%%%%%%
\section{Double branched covers and handlebody decompositions}\label{sec:dbc}

In this section we describe the double  cover of the $n$-ball $B^n$ branched along a smoothly and properly embedded codimension-two submanifold $F$. We assume that the radial distance function on $B^n$ restricts to a Morse function on the branch locus $F$. Recall the branched cover of an $n$-ball with branch locus an unknotted properly-embedded codimension-two disk is again a copy of $B^n$. By considering the gluing of this branched cover ball we show that the induced handle decomposition of $F$ determines a handle decomposition of the branched cover.

A brief description of our method is as follows: we describe the change in the branched double cover resulting from the addition of a single handle to the branch locus.
We use an imaginary ice cream scoop to remove a neighborhood of the handle  from the ball.  
Taking the double cover of a small scooped-out ball containing a $k$-handle of the branch locus results in a $(k+1)$-handle to attach to the previously constructed double branched cover.
Our main interest is in dimension 4, which we consider in \Cref{sec:dbc4}, but we begin here with a consideration of the general case, followed by a warm-up in dimension 3 in \Cref{sec:dbc3}.  Working from a suitable projection of the branch locus to $\partial B^n$ we produce either a Heegaard diagram of the double branched cover if $n=3$, or a Kirby diagram if $n=4$.

Other sources dealing with Heegaard diagrams of branched covers include \cite{josh,eli,IPY,ciprian,WB}.  Our Kirby calculus description in dimension 4 generalises those in \cite{a2016,ak,gs}, and will be used to prove that the disoriented homology of a slice surface $F$ is isomorphic to the homology of the double branched cover of $B^4$ with branch set $F$.

%%%%
\subsection{Handles and double branched covers}
%In this subsection we consider the general problem of obtaining a handle decomposition of the double branched cover of a ball given a handle decomposition of the branch locus induced by the radial distance function on the ball.  

Recall that a $k$-handle $H$ of an $n$-dimensional manifold $M$ is the image of the product $B^k\times B^{n-k}$ under an embedding $\vfi$. The attaching region of $H$ is $\vfi(\partial B^k\times B^{n-k})$ and its attaching sphere is $\Sigma:=\vfi(\partial B^k\times\{0\})$. The framing of $\Sigma$ is given by the product structure on the normal disk bundle of $\Sigma$ determined by $\vfi$.  The remainder of the boundary of $H$, $\vfi(B^k\times \partial B^{n-k})$, is its coattaching region and $\vfi(\{0\}\times \partial B^{n-k})$ is its coattaching sphere, also commonly referred to as its belt sphere.

\begin{figure}[htbp]
\centering
\includegraphics[scale=1]{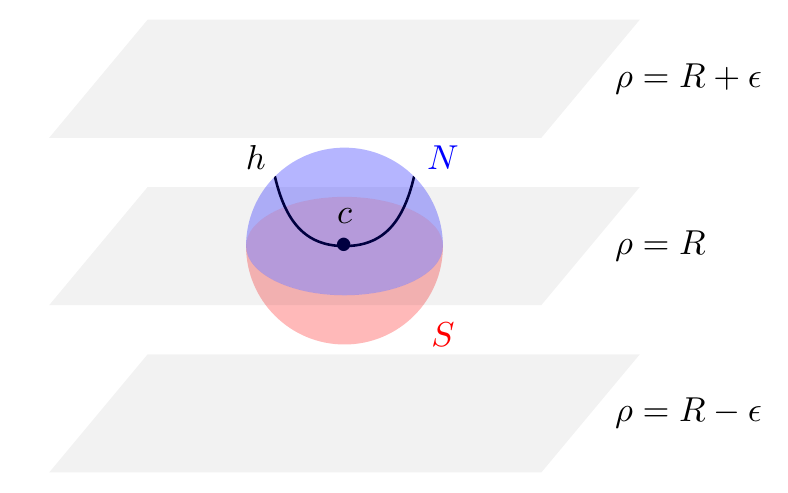}
\caption
{{\bf A ball neighborhood $D$ of a critical point.} The part of $F$ contained in $D$ is a handle $h$ corresponding to the critical point $c$. We may imagine that $D$ is attached to the sublevel set $B^n_\Rmv$ along its southern hemisphere $S$ and that its northern hemisphere $N$ is contained in $B^n_\Rpv$ by flowing the rest of these level sets into $B^n_R$.}
\label{fig:nbdD}
\end{figure}

Denote by $\rho$ the radial distance function on $B^n$. For any subset $X \subseteq B^n$ and any $r_1 < r_2$ in $[0,1]$ let $X_{r_1,r_2}$ denote $X \cap \rho^{-1}([r_1,r_2])$ and for $r \in (0,1]$ let $X_r=X_{0,r}$.
Assume that $F \subset B^n$ is a properly embedded compact codimension-two submanifold such that the restriction $\rho_F$ of $\rho$ to $F$ is Morse. Let $R$ be a critical level of $\rho_F$ that contains a single critical point $c$ whose index is $k$. Let $\veps>0$ be small enough so that $c$ is the only critical point of $\rho_F$ in $F_{\Rmv,\Rpv}$. We may choose a closed ball neighborhood $D \subset B^n_{\Rmv,\Rpv}$ about $c$ so that $h=D \cap F$ has the structure of a $k$-handle of $F$ corresponding to $c$ (see \Cref{fig:nbdD}). 
Denote the southern hemisphere $(\partial D)_R$ of $D$ by $S$ and the northern $(\partial D)_{R,1}$ by $N$. Let $C_S$ (respectively, $C_N$) be the radial  projection of the core (resp., cocore) of $h$ to $S$ (resp., $N$). The projection $h_S$ of $h$ to $S$ determines a framing $\calf_h$ of $C_S$ in $S$ as follows: the product structure on $h_S$ given by the framing of $h$ along with the normal direction to $h_S \subset S$ determine the product structure of the normal bundle of $C_S$ in $S$. Note that this framing is uniquely determined by the framing of $h$. To simplify notation we identify $S$ and $N$ with their corresponding subsets of $\partial B^n_\Rmv$ and $\partial B^n_\Rpv$ in the rest of this section.

We denote the double branched covering projection (and its restriction to any subset) by $\pi : \Sigma_2(B^n,F) \to B^n$ and the preimage of any subset $X \subseteq B^n$ under $\pi$ by $\widetilde X$.

The following theorem is the key technical result of this section.

\begin{theorem} 
\label{thm:handledbc}
With notation as above there is the following identification of double branched covering spaces:
\begin{equation}\label{Eq:dbc}
\Sigma_2(B^n_\Rpv,F_\Rpv)\cong\Sigma_2(B^n_\Rmv,F_\Rmv)\cup H,
\end{equation}
where $H$ is a $(k+1)$-handle corresponding to $\Sigma_2(D,h)$. The attaching region of $H$ in $\partial \Sigma_2(B^n_\Rmv,F_\Rmv)$ is $\widetilde S$, the preimage of $S$ under $\pi$. The attaching sphere of $H$ is $\widetilde C_S$ and its framing $\calf_H$ is given by the preimage under $\pi$ of the framing $\calf_h$.

Using identification \eqref{Eq:dbc}, the restriction of $\pi$ to $\Sigma_2(B^n_\Rpv,F_\Rpv)$ agrees with that on $\Sigma_2(B^n_\Rmv,F_\Rmv)$ away from $H$.  There are  identifications of $H$ and $D$ with $B^n$, such that the branch set $h$ corresponds to $B^{n-2} \times \{(0,0)\}$ and $\pi$ is the product of the identity on $B^{n-2}$ and the standard branched double covering projection on the normal $2$-disks.  The coattaching sphere $\widetilde C_N$ of $H$ then corresponds to $\{0^k\} \times S^{n-2-k} \times \{0\}$ and the  coattaching region $\widetilde N$ of $H$ is a regular neighborhood of $\widetilde C_N$ in $S^{n-1}$ that is diffeomorphic to $B^k \times S^{n-2-k} \times B^1$.
\end{theorem}

\begin{proof}
A standard Morse theory argument shows that $(B^n_\Rpv,F_\Rpv)\cong  (B^n_R,F_R) \cup  (D,h)$ and that $(B^n_\Rmv,F_\Rmv)\cong \overline{(B^n_R,F_R) \smallsetminus (D,h)}$ modulo corners along the equator $S \cap N$ of $D$ (see \Cref{fig:nbdD}). Here and later we suppress standard details regarding smoothing of corners. The equality of the branched covering spaces then follows from this after recognizing the branched double cover $H$ of $(D,h)$ as a $(k+1)$-handle which is the goal of the rest of the proof.

We start by choosing a convenient model for the pair $(D,h)$. Identifying $D$ with $B^n$, where the equator of $D$ is identified with the equator of $B^n$,  the handle $h$ may be identified with a part of the graph of the standard index $k$ Morse function $f : \rr^k \times \rr^{n-2-k} \times \{0\} \to \rr$, $(x,y,0) \mapsto -||x||^2+||y||^2$  (see the left side of \Cref{fig:std-model}). 
\begin{figure}[htbp]
\centering
\includegraphics[scale=1.3]{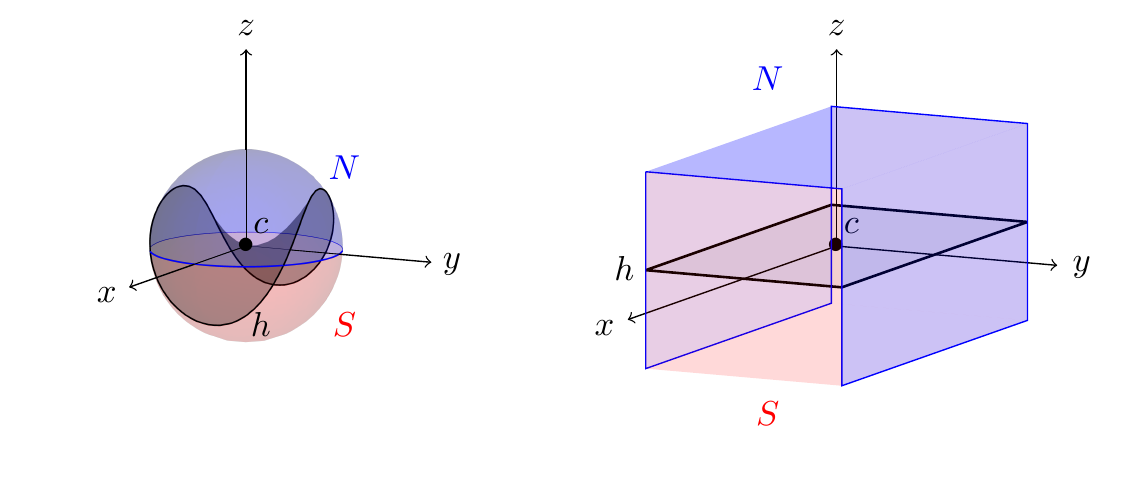}
\caption
{{\bf A standard model for the pair $(D,h)$.} Both figures show only the slice $t=0$. The left figure gives a model using a standard Morse function description of $h$ inside  $B^n$. On the right the handle has been moved to the level $z=0$ of the product ball $B$ and the subsets of $\partial B$ corresponding to $S$ and $N$ were adjusted accordingly; for $t \in (-1,1)$ the same picture describes the intersections of $S$ and $N$ with the $t$-slice.}
\label{fig:std-model}
\end{figure}
Here and below the factor $\rr^k$ gives the direction of the core of $h$ and we will denote the coordinate in this factor by $x$, the factor $\rr^{n-2-k}$ gives the direction of the cocore of $h$ and we will denote the coordinate in this factor by $y$,  the normal direction to $\rr^{n-2}$ in the domain of $f$ corresponds to the normal direction to the radial projection of $h$ into a level sphere and we will denote the coordinate in this factor by $t$, and finally the codomain of $f$ corresponds to the radial direction which we will denote by $z$. 

Applying a diffeomorphism of $B^n$ we may assume that $h$ is identified with $B^{n-2} \times \{(0,0)\}$. As a final modification we replace $B^n$ by $B:=B^k \times B^{n-2-k} \times B^1 \times B^1$ (preserving the product structure in the ambient space), where the core of $h$ corresponds to $B^k  \times \{(0^{n-2-k},0,0)\}$ and its cocore to $ \{0^k\}  \times  B^{n-2-k}  \times   \{(0,0)\}$ (see the right side of \Cref{fig:std-model}). Thus we identify $D$ with the product of $h$ with a 2-disk; the first factor of the 2-disk corresponds to the normal direction to the radial projection of $h$ and the second to the radial direction. This already shows that $H$, the double branched cover of $(D,h)$, is also diffeomorphic to $B^k \times B^{n-2-k} \times B^1 \times B^1$, with the branched covering projection $\pi$ acting nontrivially on the 2-disk $B^1 \times B^1$ given by the last two factors. This projection is essentially described by identifying $B^1 \times B^1$ with the round disk $B^2$ and using the standard branched double covering projection on that space. More precisely, we let $\pi$ be the cone (with vertex at the origin) of the orientation preserving map $\partial (B^1 \times B^1) \to \partial (B^1 \times B^1)$ that maps each of the two vertical sides $\{\pm 1\} \times [-1,1]$ diffeomorphically onto the union of the top and the right sides and each of the two horizontal sides $[-1,1] \times \{\pm 1\}$ diffeomorphically onto the union of the bottom and left sides (see \Cref{fig:disk-dbc}).
\begin{figure}[htbp]
\centering
\includegraphics{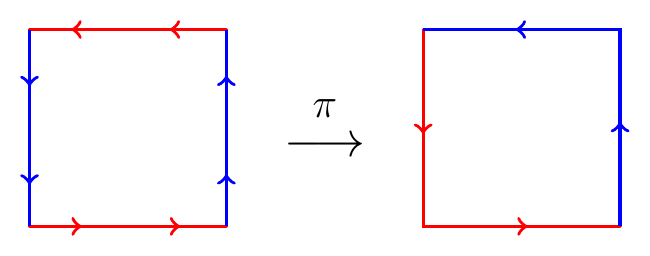}
\caption
{{\bf The double branched covering projection on $B^1 \times B^1$.} The left and right sides of the square map onto the union of right and top, and the bottom and top sides map onto the union of left and bottom, respecting orientations.}
\label{fig:disk-dbc}
\end{figure}

After the last modification we may assume that  (see the right side of \Cref{fig:std-model}):
\begin{itemize} 
\item $S$ is the union of $S_{-1}:=B^k \times B^{n-2-k} \times \big( \{-1\} \times B^1 \cup  B^1 \times \{-1\}\big)$ and $S_0:=\partial B^k \times B^{n-2-k} \times B^1 \times B^1$;

\item  $N$ is the closure of the complement of $S$ in $\partial B$, hence it is the union of $N_1:=B^k \times B^{n-2-k} \times \big( \{1\} \times B^1 \cup  B^1 \times \{1\} \big)$ and $N_0:=B^k \times \partial B^{n-2-k} \times B^1 \times B^1$.
\end{itemize}
Note that we made a choice to include one of $B^k \times B^{n-2-k} \times \{\pm 1\} \times B^1$ into $N$ and one into $S$.

Since $S_{-1}$ does not intersect the branch set, $\widetilde S_{-1}$ consists of two copies of this set, which are identified with  $B^k \times B^{n-2-k} \times B^1 \times \{\pm 1\}$. However, $S_0$ intersects the branch set in the attaching region $\partial B^k \times B^{n-2-k} \times \{(0,0)\}$ of $h$ and so $\widetilde S_0$ may be identified with $\partial B^k \times B^{n-2-k} \times B^1 \times B^1$ where $\pi$ is nontrivial on the 2-disk $B^1 \times B^1$. Hence $\widetilde S=\widetilde{S}_{-1}\cup\widetilde{S}_0$ is identified with $\partial B^k \times  B^{n-2-k} \times B^1  \times B^1 \cup B^k \times  B^{n-2-k} \times B^1 \times  \partial  B^1 \cong S^k \times B^{n-1-k}$. 
This implies that the attaching sphere of $H$ is $\partial (B^k  \times   \{(0^{n-2-k},0)\} \times B^1) $ which corresponds to $\widetilde C_S$ and its framing is given by the pull-back of the product structure on the projection of $h$ to $S$ along with the direction normal to this projection in $S$.

A similar argument as above shows that the coattaching sphere of $H$ is $\{0^k\} \times \partial  (B^{n-2-k} \times B^1) \times \{0\}=\widetilde C_N$ and its coattaching region is $\widetilde N = B^k \times \partial  (B^{n-2-k} \times B^1) \times B^1$. The restriction of $\pi$ to $ (B^k \times \partial  B^{n-2-k}) \times (B^1 \times B^1)$ is the product of the identity on $B^k \times \partial  B^{n-2-k}$ and the standard branched covering projection on the 2-disk $B^1 \times B^1$, and its restriction to $B^k \times B^{n-2-k} \times \{\pm 1\} \times B^1$ stretches each of the vertical sides of the 2-disk $B^1 \times B^1$ to the union of its right and top sides as described above. This description agrees with the one in the statement of the lemma after replacing the product ball $B$ by the round ball $B^n$.
\end{proof}

We give a more explicit description of the branched covering projection $\pi$ on the coattaching region of a handle as described in the proof of the previous lemma. This is important for understanding gluings of handles in the branched double cover of the ball of index greater than 1 as parts of their attaching regions go over coattaching regions of lower index handles.

\begin{corollary}
\label{cor:coattachingdbc}
Consider a handle $h$ of $F$ and its corresponding handle $H$ in the branched double cover $\Sigma_2(B^4,F)$ as in \Cref{thm:handledbc}. Identify $N$ with $B^k \times 2B^{n-1-k}$ where the radial projection $C_N$ of the cocore of $h$ corresponds to $\{0^k\} \times B^{n-2-k} \times \{0\}$ and the remaining direction in $2B^{n-1-k}$ is normal to the radial projection of $h$ in $N$.  Let $\Delta:=B^{n-2-k} \times B^1$ be obtained by cutting $2B^{n-1-k}$ along the annulus $(2B^{n-2-k} \smallsetminus B^{n-2-k}) \times \{0\}$; the lateral boundary $\partial_1 \Delta:= \partial B^{n-2-k} \times B^1$ of $\Delta$ corresponds to the cut (see \Cref{fig:coattach}). Then the coattaching region $\widetilde N \cong B^k \times S^{n-2-k} \times B^1$ of $H$ may be obtained from $B^k \times \Delta_\pm$ by gluing pairs of points $(x,y,z)_- \sim (x,y,-z)_+$ in  $B^k \times \partial_1\Delta_\pm$.
\end{corollary}

\begin{proof}
In the proof of the previous lemma we identified the coattaching region $\widetilde N$ of $H$ with $B^k \times \partial(B^{n-2-k} \times B^1) \times B^1$, where the coattaching sphere $\widetilde C_N$ is $\{0^k\} \times \partial(B^{n-2-k} \times B^1) \times \{0\}$. Recall that the covering transformation acts by the identity on $B^k \times B^{n-2-k}$ and by the half-turn rotation on the disk $B^1 \times B^1$. A fundamental domain for this action on the coattaching region is $B^k \times \big(B^{n-2-k} \times \{1\} \cup \partial  B^{n-2-k} \times [0,1]\big) \times B^1$, as shown in \Cref{fig:coattach}. 

\begin{figure}[htbp]
\centering
\includegraphics[scale=0.8]{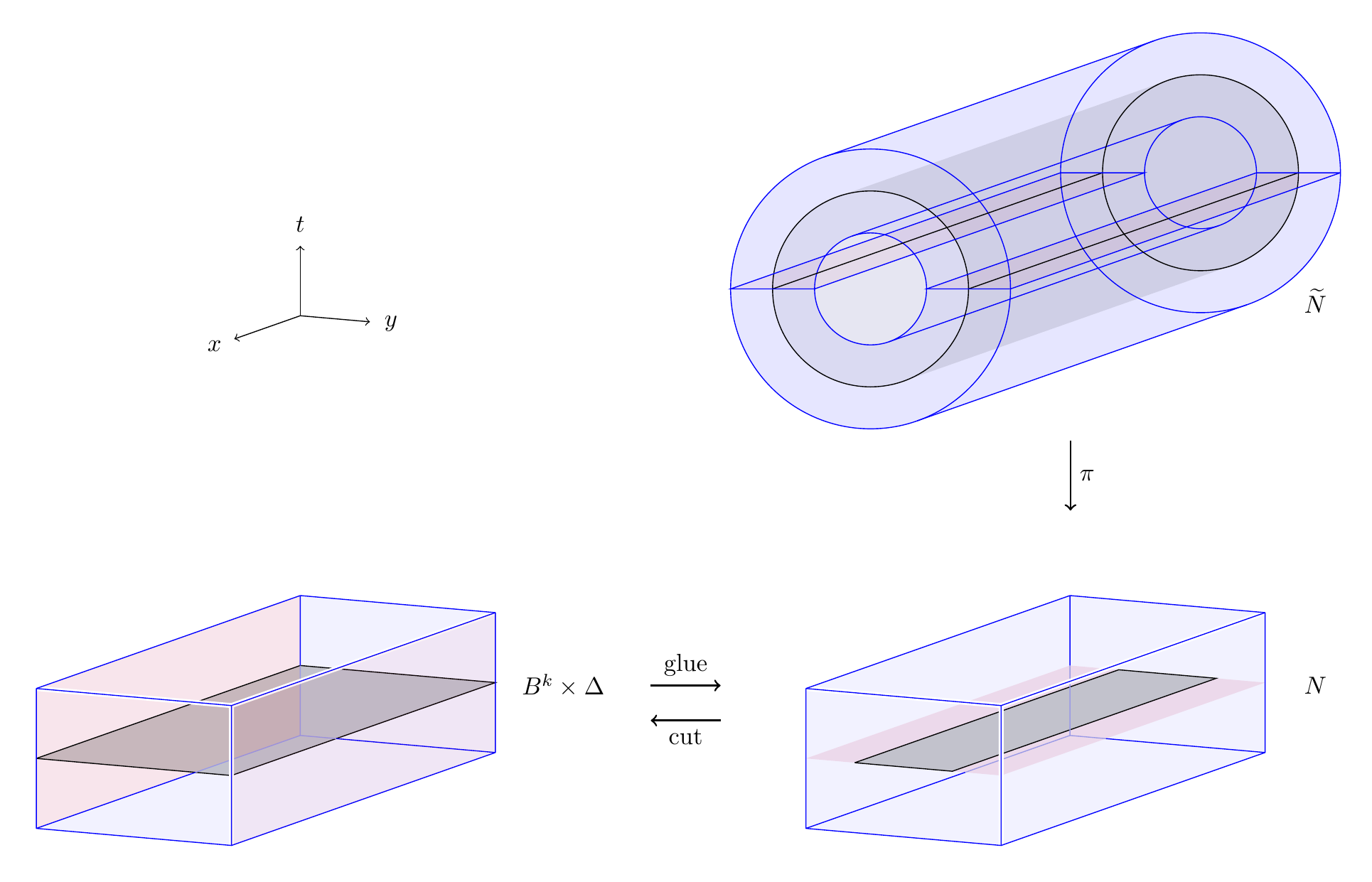}
\caption
{{\bf The coattaching region of $H$ for $n=4$ and $k=1$.} The top  picture represents a round model of $\widetilde N$ as described in \Cref{thm:handledbc}; the $z$ direction is projected into the $yt$-plane as the thickness of the annulus. 
The bottom left picture contains the fundamental domain $\Delta$ (in fact, $B^k \times \Delta$).  The bottom right picture shows $N$, and the black rectangle is the radial projection of the band (1-handle) $h$.
}
\label{fig:coattach}
\end{figure}

The branched covering projection $\pi$ maps the first set bijectively to $N_1=B^k \times \big(B^{n-2-k} \times \{1\} \times B^1 \cup B^{n-2-k} \times B^1 \times \{1\}\big)$ and the second onto $N_0=B^k \times \partial B^{n-2-k} \times B^1 \times B^1$ identifying the points $(x,y,0,z)$ and $(x,y,0,-z)$. The branch set for $\pi$ restricted to $\widetilde N$ is $B^k \times \partial B^{n-2-k} \times \{0\} \times \{0\}$. Since the maps act as identity on the first factor $B^k$, we restrict our attention to the remaining factors. We further identify the rest of $\widetilde N$ with the annulus $S^{n-2-k} \times B^1$, the fundamental domain for the action with $S^{n-2-k}_+ \times B^1 \cong B^{n-2-k} \times B^1=: \Delta$ and the branch set with the equator $S^{n-3-k} \times \{0\}$. Then $N$ (modulo $B^k$) is obtained from $\Delta$ by identifying pairs of points $(y,z) \sim (y,-z)$ in $S^{n-3-k} \times B^1$ and hence may be identified with $2B^{n-1-k}$ with the radial projection of the cocore of $h$ corresponding to $B^{n-2-k} \times \{0\}$. Conversely, $\widetilde N$ may be obtained from two copies of  $2B^{n-1-k}$ cut along the annulus $(2B^{n-2-k} \smallsetminus B^{n-2-k}) \times \{0\}$. The cut ball is diffeomorphic to $\Delta$ and the gluing of the two copies $\Delta_\pm$ identifies pairs of points $(y,z)_- \sim (y,-z)_+$ in  $(S^{n-3-k} \times B^1)_\pm$.
\end{proof}

As usual with Morse theory arguments, the assumption that there is a unique critical point in each critical level is unnecessary as the construction affects only a neighborhood of  the critical point and its preimage. In fact, we may assume all the critical points of a given index are contained in the same level set which we do in the following discussion. 
In the rest of this section we give a more detailed description of gluings of handles of small indices; we refer to the notation in the proof of \Cref{thm:handledbc}. 

%% $k=0$
\subsection{Critical points of $\rho_F$ of index $k=0$}\label{subsec:ind0}
The sublevel set $B^n_\Rmv$ is a ball that does not intersect the branch set $F$ and hence its branched double cover is the disjoint union of two $n$-balls oriented consistently with $B^n_\Rmv$, which we denote by $B^n_\pm$. Each critical point gives rise to a 1-handle connecting the two balls. Consider a 0-handle $m$ of $F$ corresponding to a critical point $c$. The attaching sphere $\{c_-, c_+\}$ of the resulting 1-handle $M$ is the preimage under $\pi$ of the radial projection of $c$ to $S$, and the attaching region $\widetilde S$ consists of two copies of $S$. Recall that $S$ is an $(n-1)$-ball centered at the radial projection of $c$; more precisely, we identify it with $B^{n-2} \times B^1$ (where the radial projection $m_S$ of $m$ to $S$ is contained in the interior of $B^{n-2}$) and the attaching map is given by
$$\vfi : (B^{n-2} \times B^1) \times  \partial B^1  \to (B^{n-2} \times B^1)_- \sqcup (B^{n-2} \times B^1)_+ \subset  B^n_- \sqcup B^n_+ ,$$
$$(y,t,z) \mapsto (y,zt)_{\sign z}.$$
Note that this is an orientable gluing. 

Alternatively, the addition of the 1-handle $M$ may be realized by gluing the balls $B^n_\pm$ along the attaching regions $ (B^{n-2} \times B^1)_\pm$ via the map
$$(y,t) \mapsto (y,-t).$$
This identifies the cocore $B^{n-2} \times B^1  \times \{0\}$ of the handle with the attaching regions and pushes one half of the handle into each of the $n$-balls. In this case it is convenient to replace the product ball $B^{n-2} \times B^1$ with the round ball $B^{n-1}$. This ball is split in half by $B^{n-2} \times \{0\}$ which we identify with the $\pi$-preimage $\widetilde m$ of the 0-handle $m$ of $F$, and we identify the boundary $\partial B^{n-1}$ with the coattaching sphere $\widetilde C_N$. For each point in $m$, which is identified with a point $y \in B^{n-2}$, we identify its corresponding points in $\widetilde C_N$ with the points $(y,\pm t)_- = (y,\mp t)_+$. When considering attachments of higher index handles we can therefore imagine that the ball $B^{n-1}$ is being inflated from the flat $B^{n-2}$, pushing the rest of the radial projection of $F$ (cut along the interior of this $B^{n-2}$) away while keeping the $y$-coordinates of the (doubled) points on the boundary fixed.

%% $k=1$
\subsection{Critical points of $\rho_F$ of index $k=1$}\label{subsec:ind1}
The sublevel set $B^n_\Rmv$ is a ball that intersects the branch locus $F$ in its 0-handles $m_i$ and hence the branched double cover of $B^n_\Rmv$ is the disjoint union $B^n_- \sqcup B^n_+$ along with a 1-handle $M_i$ connecting the two $n$-balls for each $i$. We choose to replace all the 1-handles by gluings as described above. Denote by $P$ the radial projection of $F$ into the boundary sphere $S^{n-1}$; we refer to the projections of the handles of $F$ into $P$ as handles of $P$. We assume that the 1-handles of $P$ are pairwise disjoint and that the cores of the 1-handles intersect the interiors of the 0-handles transversely in $P$. More precisely, there are two types of intersections:
\begin{itemize}
\item the attaching spheres of 1-handles lie in the union of the boundaries of 0-handles and we assume that $P$ is smooth along the attaching regions of 1-handles;

\item all other intersections are transverse and are interior to the cores of the 1-handles and to the 0-handles of $P$.
\end{itemize}
This means that the union of 0- and 1-handles of $P$ is smoothly embedded with the exception of ribbon singularities at which the cores of the 1-handles intersect the 0-handles transversely. Consider a 1-handle $h$ of $F$. The attaching circle $\widetilde C_S$ of the 2-handle $H$ corresponding to $h$ consists of two copies of $C_S$ cut along the interiors of the 0-handles of $F$ projected into $S$ as described above. If $h$ is attached to $m_i$, then  $\widetilde C_S$ intersects the coattaching sphere of $M_i$ transversely once and hence goes over $M_i$ once. If the radial projection of a point in $m_i$ (identified with $y \in B^{n-2}$) in $P$ belongs to the radial projection of the core of $h$, then  $\widetilde C_S$ intersects the coattaching sphere of $M_i$ twice (in points $(y,\pm t)_-$) and hence it goes over $M_i$ twice.

The gluing of $H=B^1 \times B^{n-3} \times B^1 \times B^1$ is determined by the attaching circle $\widetilde C_S$ corresponding to $\partial (B^1 \times \{(0,0)\} \times B^1)$ and by its framing. For $n=3$ there is a unique framing, and for $n=4$ the framing is uniquely determined by the framing of $h$, given by a parallel to the core of $h$, so by a boundary component of the $\pi$-preimage of $h_S$, the radial projection of $h$ to $S$.

For $n=4$ we recall the description of the coattaching region $\widetilde N$ of $H$ from \Cref{cor:coattachingdbc}. Choose a disk $2B^2$ in $S^3$ that intersects the radial projection of $h$ in its cocore transversely and contains this cocore in its interior as $B^1\times \{0\}$. Thicken this disk to a 3-ball $N=B^1 \times 2B^2$, where the $B^1$ factor corresponds to the core of the handle, and cut it along $B^1 \times (2B^1 \smallsetminus B^1) \times \{0\}$ to obtain $B^1 \times B^1 \times B^1$ (see \Cref{fig:coattach}). Then $\widetilde N$ is a solid torus $B^1 \times S^1 \times B^1$ obtained from two copies $(B^1 \times B^1 \times B^1)_\pm$ of the cut-up $N$ by gluing pairs of points $(x,y,z)_- \sim (x,y,-z)_+$ for $y\in \partial B^1$. If the radial projection of a 2-handle of $F$ intersects $N$ in a subset $K$, then a part of the attaching sphere of the corresponding 3-handle intersects $\widetilde N$ in two copies of $K$ cut as $N$ by the 0-handles of $F$ and glued as described above.

%%%%%%%%%%%%%%%%%%%%%%%%%%%%%%%%%%%%%%%%%%%%%%%%%%%%%%%%%%%%%%%%%%%%%%%%%%
\section{Double branched covers of the 3-ball and 3-sphere}\label{sec:dbc3}

Let $L$ be a properly embedded compact 1-manifold in the 3-ball, i.e., a tangle or a link, to which the radial distance function $\rho$ restricts to be Morse, giving a handle decomposition of $L$.  This is known as a bridge decomposition of $L$. We assume that the radial projection $P \subset S^2$ of $L$ has only ordinary double points. The bridge decomposition of $L$ induces a bridge decomposition of $P$ which then carries the same information as a diagram of $L$; we refer to double points of $P$ as crossings. In this context $0$-handles and $1$-handles are called underbridges and overbridges respectively.   We further assume that 
\begin{itemize}
\item minima of $L$ have $\rho\in(0,1/2)$ and maxima have $\rho\in(1/2,1)$,
\item all endpoints of $P$ are contained in underbridges, and
\item at each crossing, an overbridge crosses over an underbridge.
\end{itemize}

\begin{figure}[htbp]
\centering
\includegraphics[width=\textwidth]{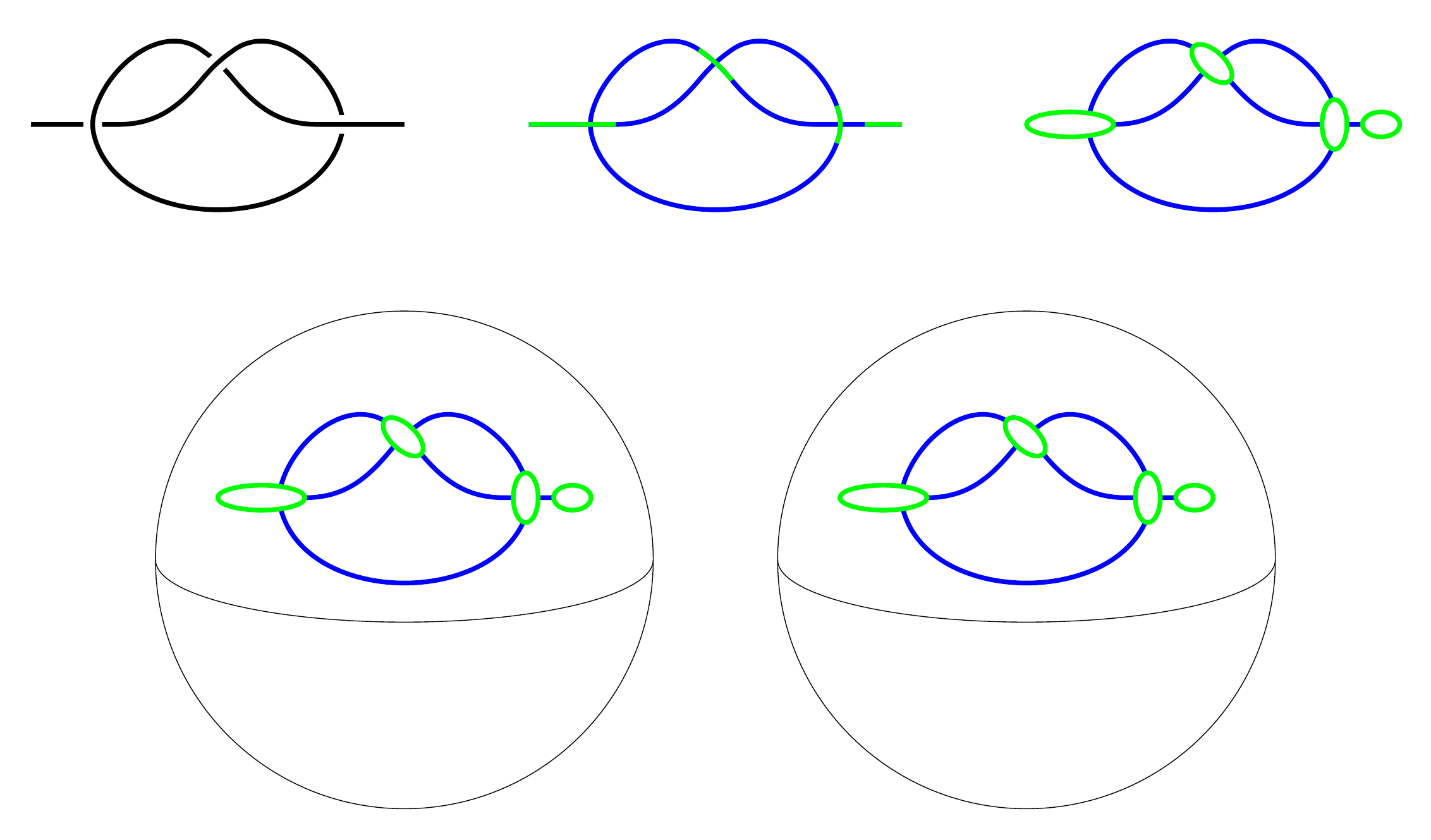}
\caption
{{\bf Double cover of a tangle in the 3-ball.} The top row shows a tangle $L$, a bridge decomposition of its projection $P$, and the associated diagram $D$ with underbridges inflated to disks.  Below these we see a handle decomposition of $\Sigma_2(B^3,L)$ with two 0-handles, four 1-handles and three 2-handles.  Matching pairs of green disks are glued preserving the direction along $L$ and reversing the normal direction.  The preimage of each blue overbridge from $D$ gives a single circle in the handlebody resulting from these disk gluings, and these are the attaching circles for the 2-handles.}
\label{fig:tangledbc}
\end{figure}

We build a handle decomposition of $\Sigma_2(B^3,L)$ using \Cref{thm:handledbc}.  We begin with a description of this which takes as a starting point any projection $P \subset S^2$ of $L$ with a chosen bridge decomposition as above.  An example is shown in Figure \ref{fig:tangledbc}.
Inflate each underbridge $u$ to a closed disk $U= B^2$ containing the underbridge as the equator $B^1\times\{0\}$, cutting any overbridge which crosses over $u$, and denote by $D$ the resulting union of disks connected by overbridge segments (see Figure \ref{fig:tangledbc}).
Let $Y$ be the oriented 3-manifold obtained from the disjoint union $B^3_-\sqcup B^3_+$ of two copies of $B^3$, each with an identical copy of $D$  in the boundary, as follows:
\begin{enumerate}
\item glue each disk $U=B^2$ in $\partial B^3_-$ to the corresponding disk in $\partial B^3_+$ by the map $(y,t)\mapsto(y,-t)$, then
\item attach a 2-handle to the resulting handlebody for each overbridge, with the attaching circle being the image in the handlebody of the union of the corresponding pair of overbridges in $\partial B^3_- \sqcup \partial B^3_+$.
\end{enumerate}
Note that each intersection of an overbridge with the boundary of a disk $U$ in $D$ results in the corresponding 2-handle attaching circle passing once over the corresponding 1-handle.  In particular, the attaching circle passes once over a 1-handle for each endpoint, and twice for each crossing. The following proposition is immediate from \Cref{thm:handledbc} and the discussion in subsections \ref{subsec:ind0} and \ref{subsec:ind1}. 

\begin{proposition}
\label{prop:tangledbc}
Let $L$ be a compact 1-manifold properly embedded in $B^3$, and let $Y$ be the 3-manifold with boundary constructed as above from a bridge decomposition of a projection of $L$.  Then $Y$ is diffeomorphic to the double  cover $\Sigma_2(B^3,L)$ of $B^3$ branched along $L$.
\end{proposition}

We would like to modify the description of $\Sigma_2(B^3,L)$ from \Cref{prop:tangledbc} to obtain a Heegaard diagram for the double branched cover of the 3-sphere along $L\subset B^3\subset S^3$.  It is convenient, though not essential, to isotope the projection $P$ so that the underbridges lie along the $x$-axis in $\rr^2\subset S^2$, and we will number them $u_0, u_1,\dots,u_g$ from left to right. 
The resulting disks in the diagram $D$ are correspondingly denoted $U_0, U_1,\dots,U_g$ from left to right. We may assume no part of $D$ lies to the left of $U_0$ with the possible exception of an arc emanating from the left endpoint of $u_0$; any other arcs may be swung across the point at infinity to the other side. 

\begin{figure}[htbp]
\centering
\includegraphics[width=\textwidth]{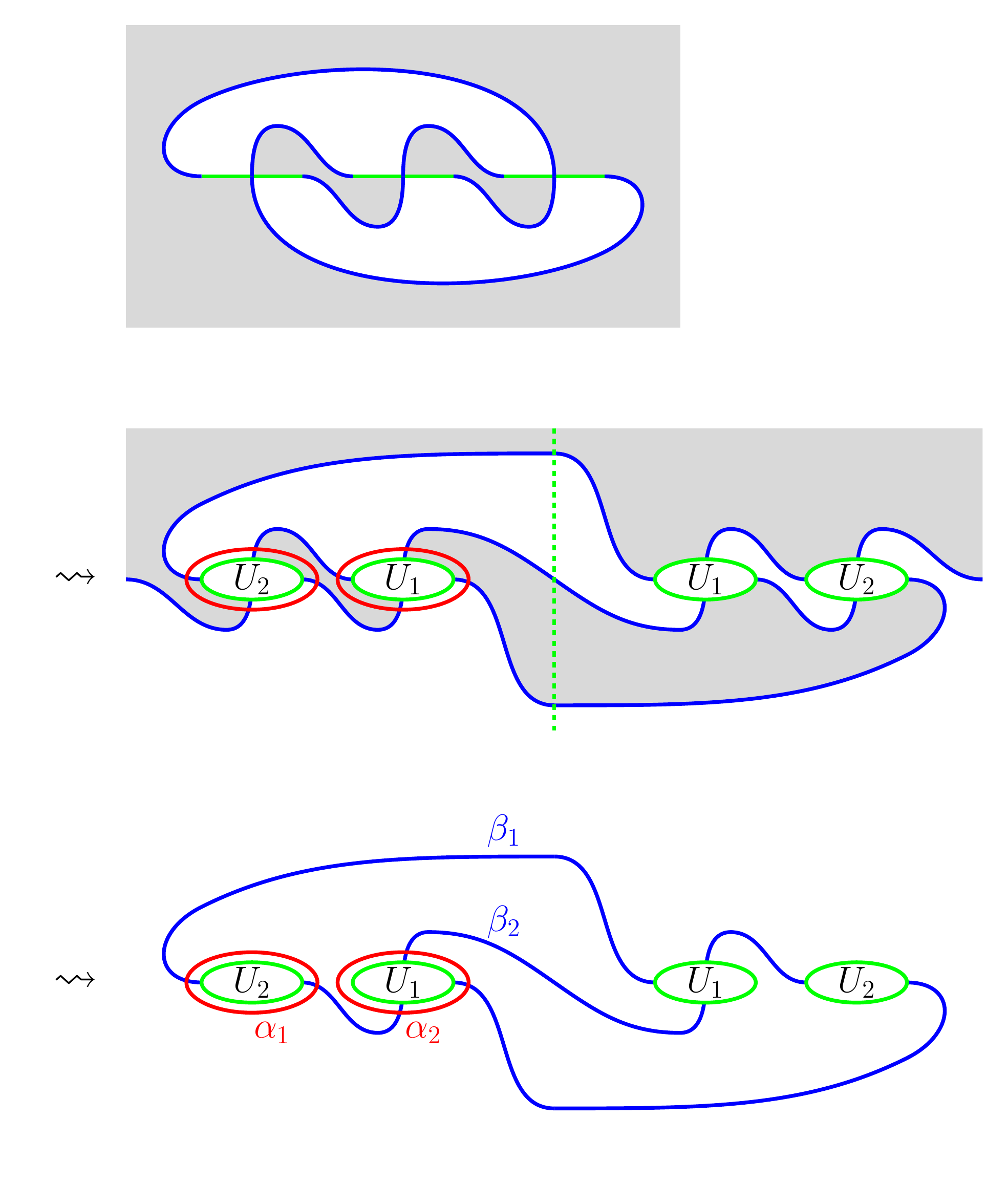}
\caption
{{\bf A Heegaard diagram for the double branched cover of the left-handed trefoil in the 3-sphere.} The chessboard colouring in the second diagram shows that the union of the blue curves is nullhomologous.}
\label{fig:HD}
\end{figure}

We now form a new planar diagram obtained as the connected sum of two copies of $D$, with the connected sum taken at the disk $U_0$, as follows. Draw one copy of $D$, with the interior of $U_0$ removed, in the right half-plane, with the disks $U_1,\dots,U_g$ drawn along the positive $x$-axis, and the boundary of $U_0$ being along the $y$-axis, with the rightmost boundary point of $U_0$ at the origin and the leftmost one at infinity.  If there is an arc emerging from this leftmost boundary point of $U_0$, redraw it as being asymptotic to the positive $x$-axis as in the second diagram of  \Cref{fig:HD}.  
For each crossing involving $u_0$, the corresponding pair of overbridge arcs should intersect the $y$-axis in a pair of points symmetric about the origin.  Draw a second copy of $D$ in the left half-plane as the rotated image of the right half-plane about the origin.  Draw a red $\alpha$ curve surrounding each of the disks in the left half-plane as in Figure \ref{fig:HD}.  The pairs of disks along the $x$-axis are now taken to be the attaching disks for 3-dimensional 1-handles, identified via reflection across the $y$-axis; thus the diagram now represents a surface $\Sigma$ of genus $g$ as the boundary of a 3-dimensional handlebody.  The blue curves coming from the overbridges form $g+1$ simple closed curves in $\Sigma$.
Let $\calh'$ denote the resulting triple consisting of the surface $\Sigma$ together with the red $\alpha$ and blue $\beta$ curves, and let $\calh$ be the triple obtained from $\calh'$  by omitting an arbitrarily chosen $\beta$ curve.

\begin{proposition}
\label{prop:heegaard}
 Let $L$ be a link in $S^3$. Then 
$\calh$ is a Heegaard diagram for the double cover of $S^3$ branched along $L$.
\end{proposition}

\begin{proof}
It is straightforward to see that $\calh'$ agrees with the description of $\Sigma_2(B^3,L)$ from \Cref{prop:tangledbc}: one first glues the two $0$-handles together using the 1-handle corresponding to $U_0$ to get a single $0$-handle.  The remaining 1-handles are indicated by the pairs of disks, and the red curves are the belt spheres of the 1-handles.  The blue curves are the attaching circles of 2-handles.  

For a link $L$, the boundary of $\Sigma_2(B^3,L)$ is a disjoint union of two 2-spheres, and the double branched cover of $S^3$ branched along $L$ can be obtained from this by attaching two 3-handles.  We claim that one of these 3-handles can be cancelled with an arbitrary choice of $\beta$ curve.   Morally, this follows from turning our construction upside-down, but 
we make a different argument.  There are $g+1$ blue $\beta$ curves on the genus $g$ surface $\Sigma$, and compressing these curves converts $\Sigma$ to a disjoint union of two spheres.  It follows that the collection of all the $\beta$ curves spans $H_1(\Sigma)$.  We claim that with an appropriate orientation, the sum of all the $\beta$ curves is nullhomologous, from which it follows that any $g$ of them span $H_1(\Sigma)$.  To see this begin by choosing a chessboard colouring of the projection $P$ of the link $L$, as in the first diagram in Figure \ref{fig:HD}.  Regions on opposite sides of an overbridge have opposite colours (shaded and unshaded).  This chessboard colouring is then inherited by the planar diagram $D$ in which the overbridges are inflated to disks.  The  Heegaard surface $\Sigma$ is obtained by taking two copies of this planar diagram, with the interiors of the inflated disks removed, and gluing them together along the boundaries of the disks.  This glues together regions from opposite sides of each underbridge, so if we choose the opposite chessboard colouring in one of the two copies 
of the planar surface being glued then the colours will match up in $\Sigma$.  Thus  the union of the $\beta$ curves bounds the union of the shaded regions.
\end{proof}

Lastly we observe that the handle decomposition of the double branched cover $Y=\Sigma_2(B^3,L)$ described in \Cref{prop:tangledbc} gives a simple way of computing the homology of $Y$ directly from a projection $P$ of $L$, equipped with a bridge decomposition. In fact, this homology is isomorphic to the disoriented homology of $L$, defined in \Cref{sec:tangle}.

\begin{proposition}
\label{prop:3d-homology}
Let $L$ be a link or tangle in $B^3$, with projection $P \subset S^2$. 
Choose a bridge decomposition of $L$ consistent with $P$ and disorientations of the overbridges, determining the data $P^\flat$. Then the homology of the disoriented chain complex $\DC_*(P^\flat)$ is isomorphic to the shifted reduced homology of $\Sigma_2(B^3,L)$, i.e.,
$$H_*(\DC_*(P^\flat)) \cong \widetilde H_{*+1}(\Sigma_2(B^3,L)).$$
\end{proposition}

\begin{proof}
This follows from the handle decomposition of $Y=\Sigma_2(B^3,L)$ described in Proposition \ref{prop:tangledbc}; we also use notation from there. Recall that 1-handles of this decomposition correspond to underbridges and 2-handles to overbridges. One of the 1-handles connects the two 0-handles $B^3_-$ and $B^3_+$, and the rest of them generate $H_1(Y)$. The relations in $H_1(Y)$ come from the 2-handles.

We label the overbridges $o_0,\dots,o_n$.  We claim that the chosen disorientation of each $o_k$ determines an orientation of the attaching circle $\beta_k$ for the corresponding 2-handle. Orient the copy of $o_k$ in $B^3_-$ consistently with $o_k$ and choose the opposite disorientation for the copy in $B^3_+$. Since for an endpoint $a$ of $o_k$ its two copies $a_\pm$ in $B^3_- \sqcup B^3_+$ are identified in $Y$, the chosen orientations match up. For a pair of endpoints $c,d$ of $o_k$ at a crossing, $c_-$ is identified with $d_+$ and $d_-$ with $c_+$, hence the chosen orientations also match up and indeed a choice of disorientation of $o_k$ determines an orientation of $\beta_k$ (see Figure \ref{fig:DHdetail}). 

\begin{figure}[htbp]
\centering
\includegraphics[scale=0.7]{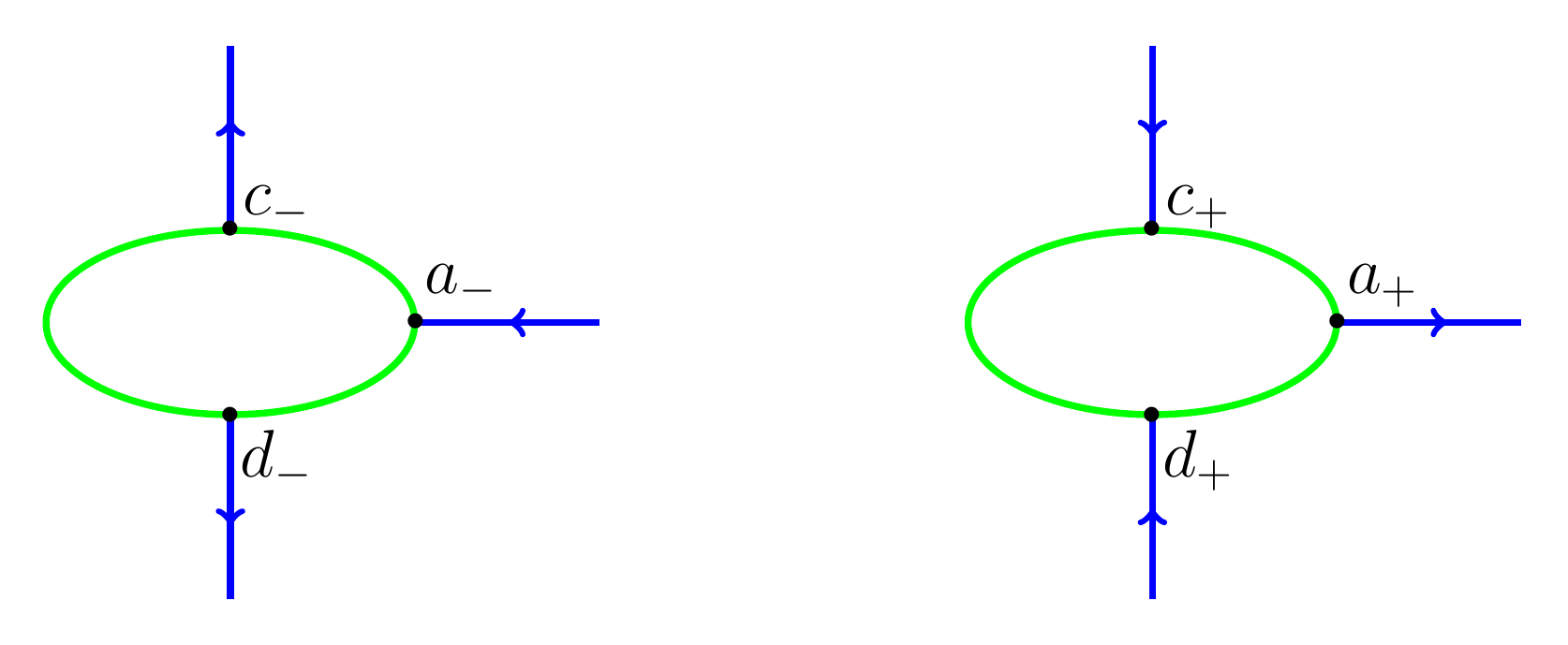}
\caption
{{\bf Disorientations of the overbridges determine orientations for the attaching circles of the 2-handles.} Arcs of overbridges $o_k$ in the two 0-handles have opposite orientations.}
\label{fig:DHdetail}
\end{figure}

We orient the 1-handles of $Y$ in such a way that the positive direction is from $B^3_-$ to $B^3_+$. Hence a 2-handle $\beta_k$ goes over a 1-handle corresponding to $u_j$ in the positive/negative direction at an endpoint $e$ of one of its subarcs if at $e$ this subarc points to/from $u_j$. 

Since $Y$ is connected and there are no 3-handles in the decomposition, the claimed isomorphism follows.
\end{proof}

Thus from \Cref{eg:trefoilDH} we see that the double cover of the 3-ball branched along a trefoil knot has first homology group isomorphic to $\zz/3\zz$, and its second homology group is a copy of $\zz$.  This is in agreement with the well-known fact that the double cover of $S^3$ branched along the left-handed trefoil is the lens space $L(3,1)$, and the double cover of $B^3$ branched along the same knot is thus obtained from $L(3,1)$ by removing two balls.

%%%%%%%%%%%%%%%%%%%%%%%%%%%%%%%%%%%%%%%%%%%%%%%%%%%%%%%%%%%%%%%%%%%%%%%%%%
\section{Double branched covers of the 4-ball}\label{sec:dbc4}

Let $F$ be a compact surface, with or without boundary, properly embedded in $B^4$.  We assume that $\rho_F$, the restriction of the radial distance function to the surface, is Morse giving a handle decomposition of $F$, and that all minima have $\rho\in(0,1/3)$, saddles have $\rho\in(1/3,2/3)$, and maxima have $\rho\in(2/3,1)$.  We further assume that
\begin{itemize}
\item the radial projection to $S^3$ restricts to an embedding on the union of $k$-handles of $F$, for each $0\le k\le2$, and we refer to the images of the handles as the handles of the projection;
\item the radial projection of the union of $0$- and $1$-handles is a ribbon-immersed surface $F_r$, and moreover all ribbon singularities are formed by 1-handles passing through 0-handles of $F_r$;
\item the radial projection $P$ of $F$ is generic, and the intersection of the interior of each 2-handle with $F_r$ is transverse.
\end{itemize}
Under these assumptions we defined a description $F_s \subset S^3$ of $F$ in \Cref{sec:surface}, which is given by the decomposition of $P$ into the ribbon-immersed surface $F_r$ and the 2-handles $d_i$. We also discussed possible singular points of $F_s$ in that section. Recall that the ribbon-immersed surface $F_r$ and the 2-handles $d_i$ may form pinch point singularities along their common boundaries. To simplify the description of the attaching spheres of 3-handles of the double branched cover we also assume that 
\begin{itemize}
\item pinch points do not occur along the boundaries of the 1-handles of $F_r$. 
\end{itemize}
Indeed, they can always be transferred along the boundary of $F_r$ by rotating the disk $d_i$ about this boundary. Hence an essential intersection of $d_i$ with some 1-handle of $F_r$ is either a component of the coattaching region of the 1-handle or is disjoint from the coattaching region, so it runs along the core of the 1-handle.

We now describe a smooth 4-dimensional handlebody $X$ diffeomorphic to the double branched cover of $(B^4,F)$, using the data above. The sublevel set $B^4_{2/3}$ is a ball that intersects the branch locus $F$ in a ribbon surface $F_{2/3}$ with projected ribbon-immersed surface $F_r$. Then the description of the branched double cover $X_2=\Sigma_2(B^4_{2/3},F_{2/3})$ is as in subsections \ref{subsec:ind0} and \ref{subsec:ind1}. Let $d$ be a 2-handle of $F$ and let $d_P$ be its radial projection into $\partial B^{4}_{2/3}$. We know from \Cref{thm:handledbc} that $d$ gives rise to a 3-handle $D$ of $X$ attached to $\partial X_2$. The attaching sphere for $D$ is the preimage of $d_P$ under the branched covering projection.
This sphere is formed by the union of the two copies of $d_P$ cut along the interiors of the 0-handles of $F_r$ in the boundaries of $B^4_\pm$, in the complement of the attaching regions of the 2-handles of $X_2$, together with the preimages of $d_P$ in the coattaching regions of the 2-handles of $X_2$. 

The coattaching region $\widetilde N_h$ of a 2-handle $H$ of $X_2$ corresponding to 1-handle $h$ of $F$ is a solid torus $B^1 \times S^1 \times B^1$ whose core circle $S^1$ is the $\pi$-preimage of the radially projected cocore of $h$. The image $N_h$ of the coattaching region under the branched covering projection $\pi$ may be identified with $B^1 \times 2B^1 \times B^1$, where the first factor corresponds to the core of $h$, the second to an extended cocore of $h$, and the third to the normal direction to $h$. With the assumptions on the projection $P$ of the surface $F$ we made above there are two types of intersections between $d_P$ and 1-handles of $F_r$. First, an arc $A$ in the boundary of $d_P$ may be glued to one of the coattaching arcs of a 1-handle $h$ of $F$. Then the component $\Delta_A$ of $d_P \cap N_h$ containing $A$ is a collar on $A$ in $d_P$; we denote the rest of the boundary of $\Delta_A$ by $A'$ (compare Figure \ref{fig:coattach} where $\Delta_A$ could be one of the two reddish rectangles in the  bottom right picture). The preimage $\pi^{-1}(\Delta_A)$ is a disk isotopic to $B^1 \times \{*\} \times B^1$ whose boundary circle is $\pi^{-1}(A')$. This disk connects the two copies $A'_\pm$ of $A'$ inside the balls $B^4_\pm$ glued along the attaching regions of the 1-handles. Since $\pi^{-1}(\Delta_A)$ intersects the core circle of the coattaching region of $H$ transversely once, the subdisk $\pi^{-1}(\Delta_A)$ of the attaching sphere $\pi^{-1}(d_P)$ goes over the 2-handle $H$ once.

The second possibility is that $d_P$ intersects $h$ in an interior arc $B$ that in $h$ runs parallel to the core of $h$. Then the component $\Delta_B$ of $d_P \cap N_h$ containing $B$ is identified with $B \times [-1,1]$ (with $B$ corresponding to $B \times \{0\}$) and $\pi^{-1}(\Delta_B)$ consists of two disks transverse to the core of the coattaching region, each capping-off one component of $(\partial \Delta_B)_\pm$.

We have now achieved the main goal of this section: a description of a handlebody corresponding to the double branched cover of a slice surface in the 4-ball.

\begin{proposition}
\label{prop:dbc4v1}
Let $F$ be a compact surface properly embedded in $B^4$ as above, and let $X$ be the 4-dimensional handlebody constructed above using the slice surface description $F_s$ of $F$.  Then $X$ is diffeomorphic to the double cover $\Sigma_2(B^4,F)$ of $B^4$ branched along $F$.
\end{proposition}

We next describe how to draw a Kirby diagram of $\Sigma_2(B^4,F)$ based on the handle decomposition from \Cref{prop:dbc4v1}.  For ribbon surfaces, this is similar to diagrams described in \cite[\S 6.3]{gs} and \cite[\S 11.3]{a2016}.  The main adjustment that needs to be made to the description above is that we need to cancel one of the two 0-handles, and draw the diagram in the boundary of the remaining 0-handle.  This is similar to what we did in the 3-dimensional case to obtain a Heegaard diagram (see \Cref{prop:heegaard}).  We begin by isotoping the radial projection $P$ of $F$ in $S^3$ to facilitate this. We assume that $P$ is contained in the upper half-space of $\rr^3 \subset S^3$ and that the 0-handles are round disks in the $xz$-plane, with their centers along a horizontal line $L$ one unit above the $x$-axis.  We then want to ``comb up" the 1- and 2-handles of $P$, so that, as much as possible, they lie above the 0-handles in the upper half space $z\ge1$ and close to the $xz$-plane.  The 1-handles (bands) are attached at their ends to the 0-handles, and pass through the 0-handles making ribbon singularities.  Away from the ends and the ribbon singularities, they are isotoped to lie close to the $xz$-plane, allowing for twisting in bands and crossings of bands over each other.  We also isotope so that the ribbon singularities all lie on the line $L$.  This gives the preferred position of the ribbon-immersed surface $F_r$. The embedded disks of the 2-handles are attached along their boundaries to the boundary of $F_r$.  Their interiors may intersect the interior of $F_r$.  Finally they may ``wrap around" the 0-handles. 
By changing the point at infinity (placing it below a chosen 0-handle and above any 2-handles wrapping around it) we may isotope $P$ so that no 2-handle wraps around a  particular 0-handle $m_0$.  We then make a further isotopy, pulling $m_0$ downwards so that it lies on the $x$-axis, below the other 0-handles.

Having isotoped the diagram in this way, we then construct the corresponding handlebody description of $X=\Sigma_2(B^4,F)$ given prior to \Cref{prop:dbc4v1}.  We inflate $m_0$ into a 3-ball in the boundary of $B^4$. Since the interior of the 3-ball becomes interior to the 4-manifold after gluing two copies of the diagram along the two copies of this 3-ball, we may consider the complement of this interior in the boundary of $B^4$, puncture the resulting boundary 2-sphere at the south pole and isotope it onto the $xy$-plane so that the boundary of $m_0$ is mapped onto the $x$-axis. This may be done without modifying the rest of the diagram which is all drawn above the $xy$-plane. This gives one copy of the diagram, corresponding to $B^4_+$. The other copy, corresponding to $B^4_-$, is obtained by revolving the first diagram about the $x$-axis so that it appears below the $xy$-plane. We have now drawn the whole diagram in a single $\rr^3$ with rotational symmetry about the $x$-axis. We inflate the remaining 0-handles of the diagram into 3-balls, remove the interiors and identify their boundaries in pairs by the reflection in the $xy$-plane. Recall that inflating cuts parts of the diagram that intersect interiors of the 0-handles.  We can push the glued 2-spheres together in the standard way to replace them by dotted circles as in \cite[Section 1.1]{a2016},\cite[Section 5.4]{gs}. The rest of the construction proceeds as described prior to \Cref{prop:dbc4v1}. The two copies of the core of each 1-handle of $P$ above and below the $xy$-plane glue to form the attaching circle for the corresponding 2-handle of $X$ whose framing is given by one component of the boundary of the annulus into which glue the two copies of the 1-handle of $P$; in fact, in the absence of 2-handles of $P$ we may consider one component of the boundary of the annulus to be the attaching circle of the 2-handle and the other to be its framing. For each 2-handle $d$ of $P$ (cut by the interiors of the 0-handles) remove from its two copies above and below the $xy$-plane the intersections with the images of coattaching regions $N_h$ for 1-handles $h$ of $P$. If a removed component $\Delta$ contains a boundary arc $A$ of $d$, the two copies $A'_\pm$ of $\overline{\partial \Delta \smallsetminus A}$ together bound a disk in the coattaching region of $h$ that goes once over the corresponding 2-handle $H$ of $X$. If a removed component $\Delta$ lies in the interior of $d$, each component of $\partial \Delta_\pm$ bounds a disk in the coattaching region of $h$ that goes once over the corresponding 2-handle $H$ of $X$.

\begin{example}[The positive unknotted real projective plane]\label{ex:proj-Kirby}
To illustrate the above results we return to the example of the unknotted real projective plane $P=\rr\pp^2$ in $B^4$ with radial projection $P_s$ given on the left of \Cref{fig:proj-plane}. Recall that we computed the disoriented homology of this surface in \Cref{ex:proj-homology}. We complete the story now by constructing a Kirby diagram for $X=\Sigma_2(B^4,P)$. 

The projection $P_s$ consists of a disk $m$ representing the 0-handle, a band with a positive half-twist representing the 1-handle $h$ that forms a ribbon singularity with the 0-handle, and a disk $d$ split into four subdisks representing the 2-handle. The red, blue and purple curves represent intersections between $d$ and $m$ and these divide $d$ into subdisks -- the faded parts of $d$ lie behind the 0-handle and the top part of the 1-handle overlaps with two subdisks. The right part of \Cref{fig:proj-plane} shows $P_s$ in the preferred position: the attachment of the 1-handle $h$ has been moved to the side and the top arc of intersection between $d$ and $m$ has been pushed away from the 1-handle. This is realized by shortening the top right subdisk and enlarging the top left subdisk of $d$ that contains the hood at the top of the figure.

\begin{figure}[htbp]
\centering
\includegraphics[scale=1.2]{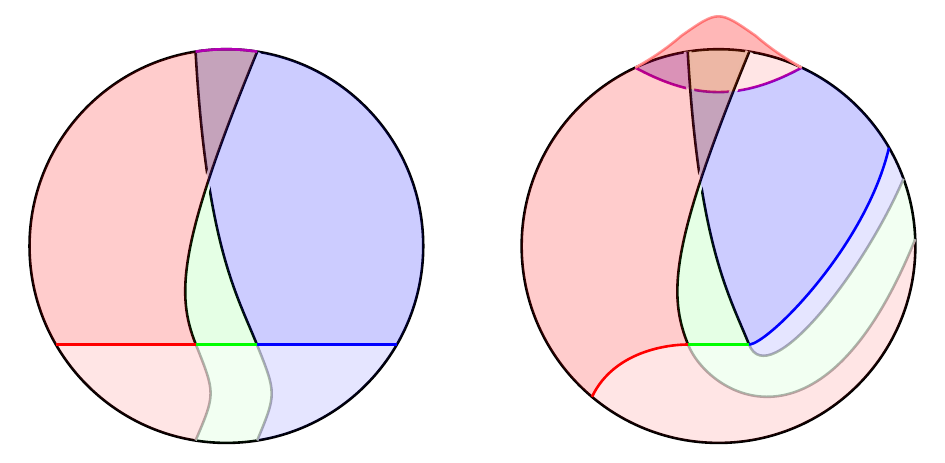}
\caption
{{\bf A radial projection of the positive unknotted real projective plane in the 4-ball.} The 0-handle (bounded by the circle) and the 1-handle (the green band) combine to give an immersed M\"{o}bius band, with a positive half-twist and a ribbon singularity (shown as a green arc).  The 2-handle consists of the red and blue disks and is split into four subdisks by its intersections with the ribbon surface shown as arcs. In the left figure, the upper two disks lie in front of the 0-handle, while the lower ones lie behind. The lower two arcs of intersection connect a pinch point with an endpoint of the ribbon singularity, and the remaining arc is bounded by pinch points. The projection on the left is not in the preferred position while the one on the right is.}
\label{fig:proj-plane}
\end{figure}

The double cover $X_2$ of $B^4$ branched along the union of the 0- and 1-handles of $P$ consists of two 0-handles $B^4_-\sqcup B^4_+$ glued along the inflated copies of the 0-handle $m$ with a single 2-handle attached. The attaching region for the 2-handle is the union of the two copies of the 1-handle $h$ of $P_s$ in the boundaries of $B^4_\pm$ cut at the ribbon singularity and pushed away by inflation of $m$. The resulting four bands form an annulus with a full positive twist. The core of the annulus is the attaching circle for the 2-handle of $X_2$ and the framing is given by either boundary component of the annulus, thus the framing coefficient is $+1$.  The Kirby diagram for $X_2$ is obtained as described above by cancelling the 1-handle with one of the 0-handles.  This results in a single 0-handle and 2-handle.  

It remains to describe the attaching sphere $S$ of the 3-handle of $X$. The part of $S$ contained in the boundary of $X_1$ away from the attching region of the 2-handle $H$ of $X$ is obtained from the two copies $d_\pm$ of $d$, cut along the 0-handle of $P_s$ and with a neighborhood of the 1-handle $h$ removed. The boundary of the resulting surface consists of two oppositely oriented framing curves of $H$ and the attaching sphere is completed by adding the two disk parallel to the core of $H$ bounded by these curves (see \Cref{fig:attach-sphere}). Of course these disks may be pushed into the coattaching region of $H$ showing that $S$ is isotopic to an unknotted 2-sphere in the boundary of the 0-handle of $X$.  We conclude that $X$ is diffeomorphic to a twice-punctured $\cc\pp^2$, as expected.
\begin{figure}[htbp]
\centering
\includegraphics[scale=1.2]{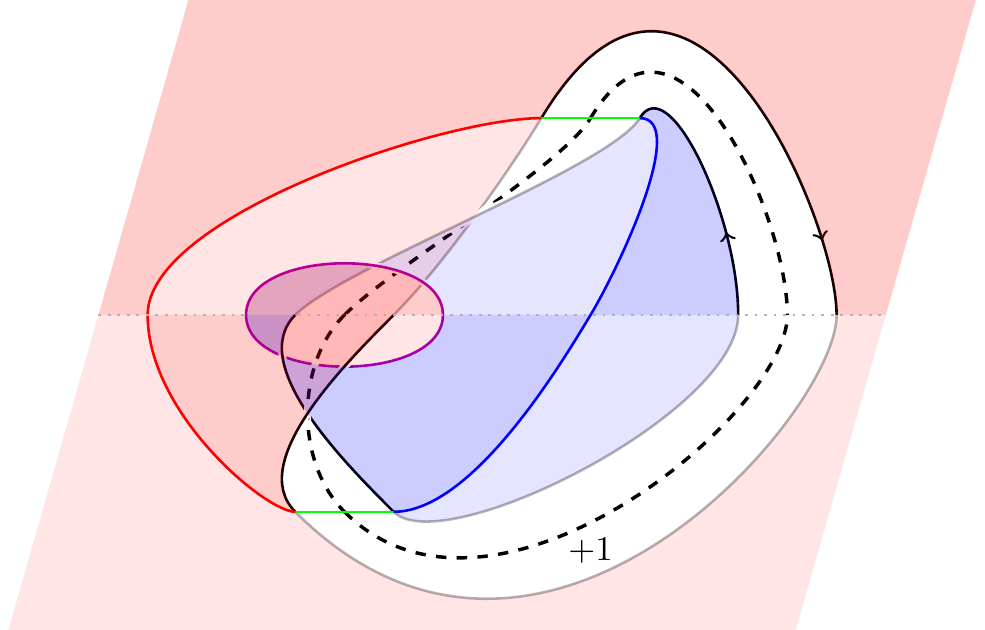}
\caption
{{\bf A Kirby diagram for $X=\Sigma_2(B^4,P)$.} The attaching circle (dashed black curve) of the 2-handle $H$ of $X$ is the $+1$-framed core of the annulus. The attaching sphere $S$ of the 3-handle is built from two copies of the subdisks of $d$ (shown with the same color scheme as in \Cref{fig:proj-plane}) with a neighborhood of the annulus removed. The solid black curves are framing curves of $H$ along which two disks parallel to the core of $H$ are attached to form $S$. The indicated orientations of the framing curves come from a choice of orientation of the visible part of $S$.}
\label{fig:attach-sphere}
\end{figure}
\end{example}

%%%%%%  Homology isomorphism
We now prove the main theorem, establishing an isomorphism between the disoriented homology of a slice surface and the homology of the 4-ball branched along the surface.
\begin{theorem}
\label{thm:4d-homology}
Let $F \subset B^4$ be a properly embedded compact surface and let $F_s \subset S^3$ be its description. 
Choose disorientations of the cores of the 1-handles and disorientations of the 2-handles of $F_s$. Then the homology of the cellular disoriented complex $\DC_*(F_s^\flat)$ is isomorphic to the shifted reduced homology of the branched double cover $\Sigma_2(B^4,F)$, i.e.,
$$H_*(\DC_*(F_s^\flat)) \cong \widetilde H_{*+1}(\Sigma_2(B^4,F)).$$
Moreover, the intersection pairing of $\Sigma_2(B^4,F)$ under this identification agrees with the GL-pairing $\lambda$ on $DH_1(F_s^\flat)$.
\end{theorem}

\begin{proof}
We first show that the homology of the double branched cover $X=\Sigma_2(B^4,F)$ with branch set a slice surface $F \subset B^4$ is isomorphic to the disoriented homology of the slice surface description $F_s \subset S^3$ of $F$. Recall that a handle decomposition of $F$ determines a handle decomposition of the projected surface $F_s$; the union of $0$- and $1$-handles of $F_s$ forms a ribbon-immersed surface $F_r$. According to \Cref{thm:handledbc} %\Cref{prop:dbc4v1} 
there is a bijection between $k$-handles of $F_s$ and $(k+1)$-handles of $X$; more precisely, the attaching sphere of a 4-dimensional handle is determined by the core of the 2-dimensional handle. Inspecting this correspondence we see that the boundary homomorphisms in the cellular disoriented complex of the surface, $\DC_*(F_s^\flat)$, and in the cellular chain complex of $X$, $C_{*+1}(X)$, agree in nonnegative dimensions. By a slight abuse we treat a handle of index $k$ as a $k$-cell. We describe below a chain equivalence inducing the claimed isomorphism. Our description relies also on the Kirby diagram described after \Cref{prop:dbc4v1}.
$$
\begin{CD}
\DC_2(F_s^\flat) @>{\partial_2^\flat}>> \DC_1(F_s^\flat)  @>{\partial_1^\flat}>> \DC_0(F_s^\flat) @>{\varepsilon}>> \zz  @>>> 0 \\
@V{\cong}V{f_2}V @V{\cong}V{f_1}V @V{\cong}V{f_0}V  @VV{f_{-1}}V @VVV \\
C_3(X) @>{\partial}>>  C_2(X) @>{\partial}>{\phantom{\partial}}>  C_1(X) @>{\partial}>>  C_0(X) @>{\varepsilon}>>  \zz
\end{CD}
$$

Recall that $X$ is built from the disjoint union of two 4-balls $X_0:=B^4_- \sqcup B^4_+$ by attaching handles. The preimage of the core of each $k$-handle of $F_s$ is a $k$-dimensional sphere in the boundary of the handlebody $X_k$, built from $X_0$ by attaching handles of index at most $k$.  Recall that the attaching sphere contains the two copies of the core in $B^4_\pm$ away from its intersection with the surface of lower index handles, connected over the coattaching regions of the corresponding handles in $X_k$.

Let $f_{-1}(1)=x_+ - x_-$, where $x_\pm$  is the generator of $C_0(X)$ corresponding to $B^4_\pm$. This makes the rightmost square commutative.

To each 0-handle $m$ of $F_s$ corresponds a 1-handle $M$ in $X$ (realized by gluing the two 4-balls in $X_0$ along the two copies of a 3-disk obtained by inflating $m$). We orient (the core of) $M$ from $B^4_-$ to $B^4_+$ and let $f_0(m)=M$, so $f_0$ sends $m$ to the oriented (core of the) 1-handle $M$. Since all the 1-handles connect the two 0-handles of $X$, the equality $\partial \circ f_0= f_{-1} \circ \varepsilon$ follows from the definition of $f_{-1}$.

For a 1-handle $h$ of $F_s$ let $c$ be its disoriented core.
At each intersection of $c$ with a 0-handle $m$ of $F_s$, the preimage $\pi^{-1}(c)=c_- \cup c_+$ of $c$ in $X_1$ goes over the corresponding 1-handle $M$ of $X$ and continues in the other ball. Orienting $c_-$ consistently with the chosen disorientation of $c$ and giving  $c_+$ the opposite disorientation, yields an oriented circle that is the attaching circle for the 2-handle $H$ of $X$ corresponding to $h$. Setting $f_1(h)=H$ it follows that $f_0 \circ \partial_1^\flat=\partial \circ f_1$ since at each point of disorientation the attaching circle goes over the 1-handle twice in the same direction.

The attaching sphere $S$ of a 3-handle $D$ of $X$ corresponding to a disoriented 2-handle $d$ of $F_s$ is obtained from the two preimages $d_\pm$ of $d$ in $X_2$. Recall that $d$ is split into faces of the graph $\Gamma=d \cap F_r$ and a disorientation of $d$ is given by a chessboard coloring of these faces. To construct $S$, change the disorientation of $d_+$ and then connect different colored faces of $\Gamma_\pm$ along the boundary of $d_\pm$ and same colored faces along the interior arcs of $\Gamma_\pm$, where all the connections are made over the handles in $X_2$. More precisely, the two copies $A_\pm$ of an arc $A$ along which $d$ is attached to a 1-handle $h$ correspond to the inclusion of a disk that goes once over the 2-handle $H$ into $S$. The sign of this contribution to the boundary is determined by the chosen orientation of $S$: if the disorientation of $d$ induces in $A$ the chosen disorientation of $h$, the sign of $H$ is positive, and negative otherwise. Similarly, the two copies $B_\pm$ of an interior arc of intersection $B \subset d$ with a 1-handle $h$ correspond to the inclusion into $S$ of two disks each of which goes once over the 2-handle $H$. The sign of this contribution may be determined from $d_-$ as before and is the same also for the other component, since both the intervening disorientations (of 1- and 2-handle) have been changed. Interpreting $S$ in $C_2(X)$ now shows that it corresponds to the disoriented 1-cycle $b^\flat$ in the definition of $\partial_2^\flat d$, proving the commutativity of the left square above.

We now turn to the pairing. Note that it is enough to establish the correspondence between pairings for ribbon surfaces. Start with a ribbon-immersed surface $F_r$ in preferred position as described in the construction of a Kirby diagram for $X$ following \Cref{prop:dbc4v1}. To simplify the discussion we additionally assume that 
no 1-handle of $F_r$ is attached along the boundary of the slice of a 0-handle bounded by the projection of a 1-handle forming a ribbon singularity (see \Cref{fig:ribbon-framing}).
Choose two disoriented cycles $a$ and $b$ in $\DC_1(F_r)$. The Gordon-Litherland pairing $\lambda_{F_r}([a],[b])$ is computed as $\lk(a,\tau b)$, where $\tau$ is the double normal push-off away from ribbon singularities and is described close to a ribbon singularity in \Cref{sec:pairing} (see \Cref{fig:pushoff}). We construct a push-off $\tau b$ that is compatible with a  ``framing'' of the curve corresponding to $b$ in $\partial X_1$, i.e., the push-off of the representative of the 2-dimensional homology class in $X$ corresponding to $b$. Recall that any disoriented 1-cycle is homologous to a linear combination of cores of 1-handles of $F_r$ whose endpoints are connected by 1-chains in the union of 0-handles. For each 1-handle $h$ of $F_r$ let its disoriented core $c_h$ be the central curve of $h$ with a chosen disorientaton. Then construct its double push-off $\tau c_h$ as follows: starting at one end of $h$ the two arcs of $\tau c_h$, one in front of $F_r$ and the other behind, both project to one side of $c_h$ in the $xz$-plane and are oriented consistently with $c_h$. These arcs can be extended along the handle by retaining relative positions of arcs with respect to $h$ as it twists and turns in space. The only exceptions to this are neighborhoods of ribbon singularities where the rule is as described in \Cref{fig:ribbon-framing}; if the projection of $\tau c_h$ along $h$ arrives to the other side of the projected core as in the picture just switch their side relative to $c_h$ by passing one over and the other under $c_h$ (alternatively one could use analogous models for the arcs on the other side). 
\begin{figure}
\includegraphics{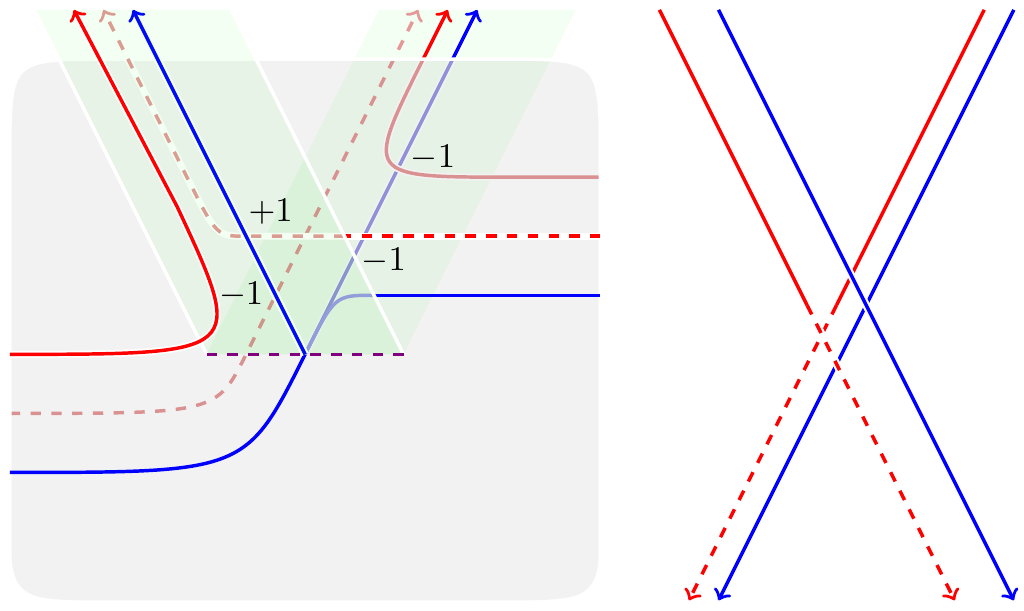}
\caption{\textbf{Comparison of pairings for a surface description in $S^3$ and for the corresponding branched double cover $X$.} The left figure shows the specific choice of framing curves (red) for a generator (blue) corresponding to a 1-handle (green) of $F_r$ near a ribbon singularity. The local contribution of the ribbon singularity to the self-pairing is $-1$. The right figure shows the corresponding attaching circle for the 2-handle (blue) and its framing (red) which also yield a local contribution of $-1$ to the linking number.}
\label{fig:ribbon-framing}
\end{figure}

To compute linking numbers we use the standard recipe of counting signs of all double points in the projection and then dividing by two (we assume all intersection points in the projection of two curves to the $xz$-plane are regular). Note first that 1-chains connecting (multiples of) disoriented cores inside the 0-handles do not contribute to the linking number $\lk(a,\tau b)$ as an intersection between $a$ and $b$ gives rise to a canceling pair of crossings between $a$ and $\tau b$. Similarly there is no contribution to $\lk(a,\tau b)$ from intersections between projections of disoriented cores and 1-chains contained in the 0-handles: if any such crossing appears, then it involves a piece of disoriented core $c_h$ pointing into/out-of a ribbon singularity, but then the same arc of a 1-chain forms an intersection also with the other piece of $c_h$ emanating from the same ribbon singularity. Since the two arcs of $c_h$ have the same orientation (pointing into/out-of the ribbon singularity) and one lies above and the other below the 0-handle, local contributions to the linking number cancel in pairs. The only contributions of arcs contained in the 0-handles of $F_r$ thus come from ribbon singularities when $[a]=[b]$. \Cref{fig:ribbon-framing} shows one possible configuration: the piece of the 1-handle $h$ lying in front of the 0-handle is on the left of the one behind. In this case the local contribution is $-1$. The other configuration is symmetric and yields local contribution $+1$.
This is consistent with the framing curve $\Phi_h$ for the attaching circle $C_h$ of the 2-handle $H$ in $X$ corresponding to $h$. The framing curve is obtained from $\tau c_h$ by keeping the ``front'' curve starting at the chosen end of $h$ in the upper half-space and rotating the other along with the diagram to the lower half-space. At a ribbon singularity this results in keeping both the front curves above the $xy$-plane and rotating the behind ones or vice versa, disregarding the parts of the curves going away from $h$ and connecting resulting arcs in the obvious way. Then any pair of crossings between $c_{h_i}$ and $\tau c_{h_j}$ corresponding to an intersection point between the projection of $c_{h_i}$ and $c_{h_j}$ results in two crossings between $C_{h_i}$ and $\Phi_{h_j}$ (see \Cref{fig:local-contribution}). 
\begin{figure}
\includegraphics[width=\textwidth]{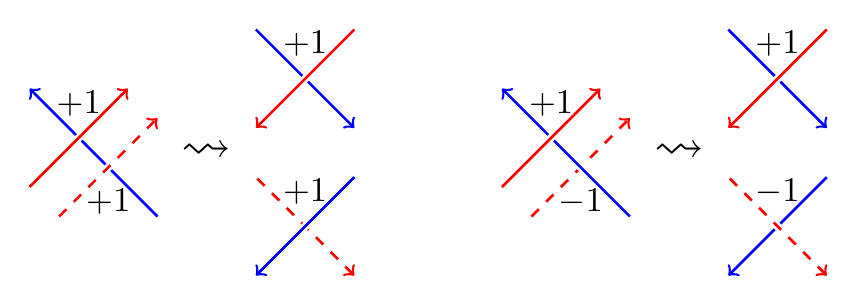}
\caption{\textbf{Comparison of pairings for a surface description in $S^3$ and for the corresponding branched double cover $X$.} The left figure corresponds to a crossing between disoriented cycles $a$ and $b$ in a surface, and the right one to a point of intersection. In the surface diagram, the blue curves represent parts of $a$ and the red ones parts of $\tau b$. In $X$, the blue curves represent parts of the attaching circle for the 2-handle corresponding to $a$ and the red ones parts of the framing curve for the 2-handle corresponding to $b$. In each case the two crossings in the surface diagram give rise to two crossings of the curves in the boundary of $X_1$ with the same local contribution to the linking number.}
\label{fig:local-contribution}
\end{figure}
The signs of these crossings agree with the signs of the original crossings since disorientations of the cores are preserved below and reversed above the $xy$-plane. The result now follows by using the above remarks and bilinearity of linking numbers. 
\end{proof}

%%%%%%%%%%%%%%%%%%%%%%%
\section{The signature formula}\label{sec:signature}

In this section we generalise Gordon and Litherland's celebrated  formula, relating the signature of a link in $S^3$ to the signature of the pairing on a spanning surface in $S^3$, to the  case of a slice surface in the 4-ball.

Let $F$ be a properly embedded surface in $S^3 \times [0,1]$ without closed components whose boundary consists of two links $\L_0 \subset S^3 \times \{0\}$ and $\L_1 \subset S^3 \times \{1\}$ (one of which could be empty). We make no assumption on orientability of $F$. We choose an orientation of the links $\L_i$ and denote the oriented links by $\vec\L_i$; recall that a link's signature is unaltered by the overall reversal of its orientation. The following proposition expresses the change in the signature of the two links in terms of the data determined by the cobordism $F$. This is a slight generalization of the signature formula in \cite{gl} and follows similarly to that. 

%Since the normal bundle of $F$ is trivial, we may choose a section $F'$ of the normal circle bundle. 
Since $F$ has the homotopy type of a 1-complex, the normal circle bundle of $F$ admits a section $F'$.
Let $\vec\L_i'$ denote the boundary links of $F'$, oriented consistently with $\vec\L_i$. Finally, let $W_F$ be the double branched cover of $S^3 \times [0,1]$ with branch set $F$.

\begin{lemma}\label{lem:sig}
With the notation as above we have
$$\sigma(\vec\L_1) - \sigma(\vec\L_0) = \sigma(W_F) + \frac 12 \left(\lk(\vec\L_0,\vec\L_0') - \lk(\vec\L_1,\vec\L_1')\right).$$
\end{lemma}

\begin{proof}
Let $\Sigma_i$ be a Seifert surface for $\L_i$. Form a (smooth) 4-sphere by adding a 4-disk to each of the boundary components of $S^3 \times [0,1]$. Then by pushing interiors of $\Sigma_i$ into the disks we may obtain a smooth surface $\widehat F$ as the union of $F$ and the pushed-in Seifert surfaces. Denoting the double branched cover of $S^4$ with branch set $\widehat F$ by $\widehat W_F$, we obtain using Novikov additivity and the $G$-signature theorem
$$\sigma(\widehat W_F)= \sigma(\vec\L_0) + \sigma(W_F) - \sigma(\vec\L_1)= - \frac 12 e(\widehat F),$$
where $e(\widehat F)$ is the normal Euler number of $\widehat F$. Recall that the normal Euler number may be computed by choosing a generic section of the normal bundle of the surface and assigning intersection numbers to intersection points by local choice of orientation of the surface and orienting the section consistently with this choice. The section $\widehat F'$ may be constructed by adding to $F'$ generic perturbations $\Sigma_i'$ of pushed-in $\Sigma_i$. As is well known, the linking number $\lk(\vec\L_0,\vec\L_0')$ is equal to the sum of local intersection numbers between $\Sigma_0$ and $\Sigma_0'$. It follows that
$$e(\widehat F)=\lk(\vec\L_0,\vec\L_0') - \lk(\vec\L_1,\vec\L_1'),$$
which proves the claimed formula.
\end{proof}

If the surface $F \subset S^3 \times [0,1]$ projects injectively to the sphere, giving an embedded cobordism between the links, the signature of the double branched cover manifold may be computed from the Gordon-Litherland pairing $\lambda_{p(F)}$. Also the links $\vec\L_i'$ may be replaced by nearby parallels of $\vec\L_i$ on the projected image of $F$.

\begin{proposition}\label{prop:sig}
Let $F \subset S^3 \times [0,1]$ be a properly embedded surface such that the restriction of the projection $p$ along the interval to $F$ is an embedding. Then for any choice of orientations of the boundary links $\vec\L_i \subset S^3 \times \{i\}$ of $F$ we have
$$\sigma(\vec\L_1) - \sigma(\vec\L_0) = \sigma(\lambda_{p(F)}) + \frac 12 \left(\lk(\vec\L_0,\vec\L_0^F) - \lk(\vec\L_1,\vec\L_1^F)\right),$$
where $\vec\L_i^F$ is a nearby parallel of $\vec\L_i$ on $p(F)$.
\end{proposition}

\begin{proof}
This follows immediately from the above lemma after noting that $\lambda_{p(F)}$ is the intersection pairing of $W_F$ and that since $F$ is a graph of $p(F)$, we can choose a section $F'$ for which $\L_i'$ is homotopic to $\L_i^F$. Indeed, a section $F'$ can be constructed starting with $\L_0'=\L_0^F$ and pushing $p(F)$ (with the collar between $\L_0$ and $\L_0^F$ removed) slightly below $F$. This has to be completed by adding a collar on the image of $\L_1$ that interpolates to $S^3 \times \{1\}$. Clearly $\L_1'$ is then homotopic to $\L_1^F$.
\end{proof}

Consider now a general slice surface $F \subset B^4$ with boundary link $\L$. We continue assuming that $F$ has no closed components. Let $F_s \subset S^3$ be a description of $F$ and $\lambda_{F_s}$ the corresponding pairing on $DH_1(F_s)$. Recall that $F_s$ consists of a ribbon surface description $F_r$ and a separated sublink of $\partial F_r$ consisting of those components that are in $F$ capped-off. Let $\L^F$ be a nearby parallel of $\L$ on $F_r$.

\begin{theorem}\label{thm:sig}
Let a link $\L$ be the boundary of a slice surface $F \subset B^4$. Then for any choice $\vec\L$ of orientation for $\L$ its signature is given by
$$\sigma(\vec\L)=\sigma(\lambda_{F_s}) - \frac 12 \lk(\vec\L,\vec\L^F),$$
where $\vec\L^F$ is oriented consistently with $\vec\L$.
\end{theorem}

\begin{proof}
We may assume the radial distance function in $B^4$ induces a Morse function on $F$ so that the ball $D_0$ of radius $1/3$ contains exactly all critical points of index $0$ and the radial shell $E_1$ between $1/3$ and $2/3$ contains exactly all critical points of index $1$. Then the part of $F$ contained in $D_1=D_0 \cup E_1$ is a ribbon surface. We further assume that the interior arcs of ribbon sigularities of $F_r$ are contained in the 0-handles.

Since only the double branched cover of $E_1$ may have nontrivial signature, it follows by Novikov additivity that the signature of the branched cover of $E_1$ equals to the signature of the branched cover of $D_1$ and to that of $B^4$, which is equal to $\sigma(\lambda_{F_r})=\sigma(\lambda_{F_s})$. The result now follows from \Cref{prop:sig} after noting that the lower boundary of the intersection of $F$ with $E_1$ is a 0-framed unlink and that 
%. That $\L'$ may be replaced by $\L^F$ follows as in the proof of the previous proposition since 
the radial projection restricted to $F \cap E_1$ is an embedding. 
\end{proof}

%%%%%%%%%%%%%%%%%%%%%%%%%%%%%%%%%%%%%%%%%%%%%%%%%%%%%%%%%%%%%%%%%%%%%%%%%%

%\clearpage

\bibliographystyle{amsplain}
\bibliography{GLslice}

\end{document}